\documentclass{amsart}
\usepackage{graphics,verbatim, amsmath, amssymb, amsthm, amsfonts}	
\input epsf         

\newtheorem{proposition}{Proposition}[section]
\newtheorem{theorem}[proposition]{Theorem}
\newtheorem{corollary}[proposition]{Corollary}
\newtheorem{lemma}[proposition]{Lemma}
\newtheorem{prop}[proposition]{Proposition}
\newtheorem{cor}[proposition]{Corollary}
\newtheorem{conj}[proposition]{Conjecture}
\newtheorem*{gVectRecurr*}{Conjectured recurrence for $\mathbf{g}$-vectors}

\theoremstyle{definition}
\newtheorem{example}[proposition]{Example}

\theoremstyle{remark}
\newtheorem{remark}[proposition]{Remark}

\numberwithin{equation}{section}

\usepackage[usenames]{color}

\newcounter{margincounter}
\newcommand{\ep}{\epsilon}
\newcommand{\cl}{\operatorname{cl}}
\newcommand{\nc}{\operatorname{nc}}
\newcommand{\NC}{\operatorname{NC}}

\newcommand{\Cat}{\operatorname{Cat}}

\newcommand{\cov}{\mathrm{cov}}
\newcommand{\covers}{{\,\,\,\cdot\!\!\!\! >\,\,}}
\newcommand{\covered}{{\,\,<\!\!\!\!\cdot\,\,\,}}
\newcommand{\set}[1]{{\left\lbrace #1 \right\rbrace}}
\newcommand{\pidown}{\pi_\downarrow}
\newcommand{\piup}{\pi^\uparrow}
\newcommand{\br}[1]{\langle #1 \rangle}
\newcommand{\A}{{\mathcal A}}
\newcommand{\EL}{{\mathcal L}}
\newcommand{\F}{{\mathcal F}}
\newcommand{\G}{{\mathcal G}}
\newcommand{\join}{\vee}
\newcommand{\meet}{\wedge}
\newcommand{\leftq}[2]{\!\!\phantom{.}^{#1} {#2}}
\newcommand{\closeleftq}[2]{\!\!\phantom{.}^{#1}\! {#2}}
\newcommand{\Pge}{{\Phi_{\ge -1}}}
\newcommand{\cm}{\parallel}
\newcommand{\ck}{^\vee}
\newcommand{\letw}{\le_{\mathrm{tw}}}
\newcommand{\Alg}{\mathrm{Alg}}
\newcommand{\toname}[1]{\stackrel{#1}{\longrightarrow}}
\newcommand{\ccc}{\gamma}
\newcommand{\groot}{\mathbf{g}_\mathrm{root}}
\newcommand{\gweight}{\mathbf{g}_\mathrm{weight}}
\newcommand{\ZZ}{\mathbb{Z}}
\newcommand{\QQ}{\mathbb{Q}}
\newcommand{\RR}{\mathbb{R}}

\newcommand{\Proj}{\mathrm{Proj}}
\newcommand{\Clust}{\mathrm{Clust}}

\newcommand{\onto}{\twoheadrightarrow}
\newcommand{\isomorph}{\cong}

\newcommand{\FFF}{\mathbb{F}}
\DeclareMathOperator{\Span}{Span}

\begin{document}
\title{Cambrian fans}

\author{Nathan Reading}
\address{Mathematics Department\\
University of Michigan\\
Ann Arbor, MI 48109-1043\\
USA}
\curraddr[Nathan Reading]{
Department of Mathematics\\
North Carolina State University\\
Harrelson Hall 225, Box 8205\\
Raleigh, NC 27695-8205\\
USA
}
\author{David E Speyer}
\curraddr[David Speyer]{
Department of Mathematics, Room 2-332 \\
Massachusetts Institute of Technology \\
77 Massachusetts Avenue \\
Cambridge, MA 02139 USA \\}
\email{nathan\_reading@ncsu.edu, speyer@math.mit.edu}


\thanks{Nathan Reading was partially supported by NSF grant DMS-0502170.
David E Speyer was supported by a research fellowship from the Clay Mathematics Institute. }

\begin{abstract}
For a finite Coxeter group~$W$ and a Coxeter element~$c$ of $W,$ the $c$-Cambrian fan is a coarsening of the fan defined by the reflecting hyperplanes of~$W\!$.
Its maximal cones are naturally indexed by the $c$-sortable elements of~$W\!$.  
The main result of this paper is that the known bijection $\cl_c$ between $c$-sortable elements and $c$-clusters induces a combinatorial isomorphism of fans.
In particular, the $c$-Cambrian fan is combinatorially isomorphic to the normal fan of the generalized associahedron for~$W\!$. 
The rays of the $c$-Cambrian fan are generated by certain vectors in the $W$-orbit of the fundamental weights, while the rays of the $c$-cluster fan are generated by certain roots.
For particular (``bipartite'') choices of~$c$, we show that the $c$-Cambrian fan is linearly isomorphic to the $c$-cluster fan.
We characterize, in terms of the combinatorics of clusters, the partial order induced, via the map $\cl_c$, on $c$-clusters by the $c$-Cambrian lattice.
We give a simple bijection from $c$-clusters to $c$-noncrossing partitions that respects the refined (Narayana) enumeration. We relate the Cambrian fan to well known objects in the theory of cluster algebras, providing a geometric context for  $\mathbf{g}$-vectors and quasi-Cartan companions.
\end{abstract}

\maketitle

\setcounter{tocdepth}{1}
\tableofcontents

\section{Introduction}
\label{intro}
Recent research in combinatorics has focused on the relationship between various objects counted by the $W$-Catalan number, $\Cat(W)$, for~$W$ a finite Coxeter group.
This number, which has a simple formula in terms of fundamental numerical invariants of $W$, has arisen separately in a wide variety of the many fields with connections to Coxeter groups.
These unexplained numerical coincidences have led to efforts to discover deeper mathematical connections between the different fields.

One set counted by $\Cat(W)$ is the set of \emph{clusters} \cite{ga,ca2} in the root system $\Phi$ associated to~$W\!$.
A cluster is a collection of roots in $\Phi$ that are ``compatible'' in a sense which will be made precise in Section~\ref{cluster sec}.
The positive linear spans of clusters are the maximal cones in a complete simplicial fan which we refer to as the \emph{cluster fan} and whose dual polytope is called the \emph{generalized associahedron} for $W\!$.
Clusters of roots get their name from cluster algebras. 
Although it is surprising \emph{a priori} that a cluster algebra should have anything to do with a Coxeter group, cluster algebras of finite type turn out to have a classification~\cite{ca2} that exactly matches the classification of finite crystallographic root systems $\Phi$.
Clusters of roots in $\Phi$ turn out to encode the combinatorics of the corresponding cluster algebra.
For a very gentle introduction to cluster algebras and to ``$W$-Catalan'' combinatorics, see~\cite{rsga}.
For a more advanced survey, see~\cite{CDM}.

Another set counted by the $W$-Catalan number \cite{Bessis,ga,MRZ,picantin,Reiner} is the set of \emph{noncrossing partitions} associated to~$W$.
Both the name and the earliest examples of noncrossing partitions come from algebraic combinatorics (see e.g.\ \cite{Kreweras,Reiner}), while both the general definition and important applications arise from geometric group theory.
Specifically, noncrossing partitions are a powerful tool in the theory of Artin groups \cite{Bessis,BWKpi}.
For an accessible introduction to this application, focusing on the special case of the symmetric group (and thus the braid group), see~\cite{McCammond}, which also discusses other applications of noncrossing partitions to free probability and combinatorics.
The definitions of both noncrossing partitions and clusters involve the choice of a Coxeter element~$c$ for $W,$ and to emphasize this fact we will refer to them as $c$-noncrossing partitions and $c$-clusters. 

A third combinatorial set counted by $\Cat(W)$ is the set of \emph{nonnesting partitions} (antichains in the root poset of $\Phi$), which will not play a role in the current paper.
These objects arose in several closely related contexts, including double affine Hecke algebras (rational Cherednik algebras), two sided cells and coinvariant rings.  See~\cite[Lecture 5]{rsga} for a gentle introduction and for references.

One of the main results of~\cite{sortable} is a bijective proof that $c$-clusters and $c$-noncrossing partitions are equinumerous.  
The proof begins with the definition of a fourth set counted by $\Cat(W)$, the set of $c$-\emph{sortable elements} of~$W$.
Bijections are then given from $c$-sortable elements to $c$-noncrossing partitions and from $c$-sortable elements to $c$-clusters.
Sortable elements and the bijections are defined simply without reference to the classification of finite Coxeter groups, but the proofs that these are bijections rest on several lemmas which are proved type by type using the classification.

Sortable elements have their origins in the lattice theory of the weak order.
Specifically, the \emph{$c$-Cambrian congruence} is a certain lattice congruence $\Theta_c$ on the weak order on~$W\!$ whose congruence classes are counted \cite{cambrian, sort_camb} by $\Cat(W)$. 
The $c$-Cambrian congruence classes are (by a general fact about congruences of a finite lattice) intervals in the weak order.
The quotient lattice $W/\Theta_c$ (the \emph{Cambrian lattice}) is isomorphic to the restriction of the weak order to the minimal elements of the $c$-Cambrian congruence classes.
These minimal elements turn out to be exactly the $c$-sortable elements~\cite{sort_camb}.

As a special case of a construction given in~\cite{con_app}, the congruence $\Theta_c$ defines a complete fan~$\F_c$, called the \emph{$c$-Cambrian fan}, whose maximal cones correspond to $c$-Cambrian congruence classes.
The fan~$\F_c$ is a coarsening of the fan defined by the reflecting hyperplanes of~$W\!$.
The goal of this paper is to understand in detail the polyhedral geometry of the bijection between the maximal cones of~$\F_c$ (indexed by $c$-sortable elements) and the maximal cones of the $c$-cluster fan (defined by $c$-clusters).
The bijection $\cl_c$ from $c$-sortable elements to $c$-clusters was defined in~\cite{sortable} without reference to the fan~$\F_c$ or the $c$-cluster fan.
The key result of this paper is the following theorem, a natural strengthening of the statement that $\cl_c$ is a bijection.

\begin{theorem}
\label{cl iso}
Let~$W$ be a finite Coxeter group and let~$c$ be a Coxeter element of~$W\!$.
Then the $c$-Cambrian fan~$\F_c$ is simplicial and the map $\cl_c$ induces a combinatorial isomorphism from~$\F_c$ to the $c$-cluster fan. 
\end{theorem}
By a combinatorial isomorphism of simplicial fans, we mean a combinatorial isomorphism of the simplicial complexes obtained by intersecting with the unit sphere. Because the fans are simplicial, this is equivalent to requiring that there be a piecewise linear homeomorphism from $\RR^n$ to itself, linear on each face of~$\F_c$, carrying the cones of~$\F_c$ to the cones of the $c$-cluster fan. The isomorphism of Theorem~\ref{cl iso} is typically only \emph{piecewise} linear.
However, we show that for any $W$, there exists a special ``bipartite'' choice of $c$ such that the $c$-Cambrian fan and the $c$-cluster fan are linearly isomorphic. 
This result (Theorem~\ref{L iso}) verifies the first statement of \cite[Conjecture~1.4]{cambrian}.

Theorem~\ref{cl iso} shows that $c$-sortable elements are not simply in bijection with $c$-clusters, but define the same underlying combinatorial structure.
This is a particularly surprising result as the Cambrian fan and cluster fan are defined in very different ways: 
the Cambrian fan is defined by removing walls of the fan defined by the reflecting hyperplanes while the cluster fan is defined by choosing certain rays, all of which are normal to reflecting hyperplanes, and specifying which rays lie in common cones.
Moreover, the Cambrian fan contributes combinatorial structure which is not present in the cluster fan, including a poset (in fact lattice) structure which interacts well with the fan structure, as well as a notion of projection to standard parabolic subgroups.

The definition of sortable elements is valid for infinite Coxeter groups.
The theory of sortable elements and Cambrian lattices/fans can be extended to infinite Coxeter groups and we will describe this program in detail in a future paper.   In particular, this future paper will provide uniform proofs, valid for all finite and infinite Coxeter groups, of the results which were proved using type by type analysis in \cite{sortable} and \cite{sort_camb}.  
Cluster algebras of infinite type are not as well understood as cluster algebras of finite type, and one of the key motivations of this paper is to establish connections between Cambrian lattices/fans which can be generalized to give new insights into cluster algebras of infinite type.
%Directly defining clusters for infinite Coxeter groups, on the other hand, appears to be very difficult because the absence of any analogue of \cite[Theorem 2.6.1]{ga} makes it unclear how compatibility ought to be defined.

Theorem~\ref{cl iso} leads to further results which we now describe. The $c$-cluster fan is combinatorially isomorphic to the normal fan of the generalized associahedron for~$W\!$.
Thus, Theorem~\ref{cl iso} implies that the $c$-Cambrian fan is combinatorially isomorphic to the normal fan of the generalized associahedron.
This isomorphism, combined with the close structural relationship which exists between the $c$-Cambrian lattice and the $c$-Cambrian fan (see Section~\ref{camb fan sec}), implies that the Hasse diagram of the $c$-Cambrian lattice is combinatorially isomorphic to the $1$-skeleton of the generalized associahedron (Corollary~\ref{Hasse}).
This confirms \cite[Conjecture~1.2.a]{cambrian}.

The $c$-Cambrian lattice induces a partial order (the {\em $c$-cluster lattice}) on $c$-clusters via the map $\cl_c$.  
We characterize this partial order in terms of the combinatorics of clusters (Theorem~\ref{cl poset iso}), generalizing and proving the second statement of \cite[Conjecture~1.4]{cambrian}.
The $c$-cluster lattice inherits many useful properties from the $c$-Cambrian lattice (see Corollary~\ref{consequences}), including the property that any linear extension of the $c$-cluster poset is a shelling of the $c$-cluster fan.
As a consequence of this shelling property we obtain, in Section~\ref{lattice sec}, a bijective proof of the fact that the $k^{\textrm{th}}$ entry of the $h$-vector of the $c$-cluster fan coincides with the number of $c$-noncrossing partitions of rank~$k$. 
(This number is called the $k^{\textrm{th}}$ Narayana number associated to~$W\!$.)

In Section~\ref{bij sec} we give a purely geometric description (Theorem~\ref{geom bij}) of a bijection between $c$-clusters and $c$-noncrossing partitions, in the case where $c$ is a bipartite Coxeter element.
This result draws on a ``twisted'' version of the $c$-cluster poset as well as the linear isomorphism, mentioned above, between the $c$-Cambrian fan and the $c$-cluster fan.
Another connection between noncrossing partitions and clusters has arisen recently.
Brady and Watt~\cite{BWlattice} construct a simplicial fan associated to $c$-noncrossing partitions (for bipartite~$c$) and extend their construction to produce the $c$-cluster fan.
Athanasiadis, Brady, McCammond and Watt~\cite{ABMW} use the construction of~\cite{BWlattice} to give a bijection between clusters and noncrossing partitions.  
Their proof uses no type by type arguments and provides a different bijective proof that the $k^{\textrm{th}}$ entry of the $h$-vector of the $c$-cluster fan coincides with the number of $c$-noncrossing partitions of rank~$k$. 
The bijection of~\cite{ABMW} incorporates elements which are similar in appearance to the constructions of the present paper (see Remark~\ref{ABMW remark}), but many details of the relation between the two theories remain unclear.

The results described above resolve all conjectures from~\cite{cambrian} except for Conjecture 1.1. This last result has been established by Hohlweg and Lange~\cite{HL} for types $A$ and $B$ and will be proven for all types in a future paper by Hohlweg, Lange and Thomas~\cite{HLT}. 

The relationship between the $c$-Cambrian fan and the $c$-cluster fan has several consequences for the theory of cluster algebras, which we describe in more detail in Section~\ref{gVectorSection}. 
When~$W$ admits a crystallographic root system~$\Phi$, Fomin and Zelevinsky associate a cluster algebra $\Alg(\Phi)$ to~$W\!$. 
This is a commutative algebra with certain specified elements, called cluster variables, and certain specified subsets of these variables, ordinarily called clusters. 
We will call these subsets of variables \emph{algebraic clusters}, to distinguish them from the combinatorially defined clusters which are certain sets of roots.   

The root system~$\Phi$ and a choice of Coxeter element~$c$ specify a certain algebraic cluster~$t_c$ of $\Alg(\Phi)$.  
There is a bijection between cluster variables and almost positive roots, such that the elements of~$t_c$ are taken to the negative simple roots and such that algebraic clusters are taken to $c$-clusters. 
To any cluster variable $x \in \Alg(\Phi)$ and any cluster~$t$, Fomin and Zelevinsky associate to the pair $(x,t)$ two vectors in~$\ZZ^n$:  the \emph{denominator vector} and the \emph{$\mathbf{g}$-vector} of~$x$ with respect to~$t$. 
It is shown in~\cite{ca2} and~\cite{ccs} that, when $t=t_c$ the denominator vector is found by expressing the corresponding root in the basis of simple roots. 
We show, in the case where $c$ is bipartite, that the $\mathbf{g}$-vector is given by expressing the corresponding ray of the $c$-Cambrian fan in the basis of fundamental weights; this result would follow for other $c$ if we knew Conjecture~7.12 of \cite{ca4}. 
We also provide a geometric context for the notion of \emph{quasi-Cartan companions}, defined in~\cite{capsm}. 
 
 In the following four sections, we lay out the necessary background concerning Coxeter groups, sortable elements and Cambrian congruences. We also give more precise statements of several of the results described in this introduction. In section~\ref{ray sec}, we begin presenting our proofs.

\section{The weak order}
\label{weak sec}
This section covers preliminary results on finite lattices and in particular on the weak order on a finite Coxeter group.
We assume that the reader is familiar with the most basic definitions of Coxeter groups and lattices.
Details about lattices are found in~\cite{Gratzer} and details about Coxeter groups are found in \cite{Bj-Br,Bourbaki,Humphreys}.

A {\em join-irreducible} element of a finite lattice~$L$ is an element which covers exactly one other element.
A {\em meet-irreducible} element of~$L$ is an element which is covered by exactly one other element.
A {\em homomorphism} from the lattice $L_1$ to the lattice $L_2$ is a map $\eta:L_1\to L_2$ with $\eta(x\meet y)=\eta(x)\meet\eta(y)$ and $\eta(x\join y)=\eta(x)\join\eta(y)$ for every $x,y\in L_1$.
The condition that $\eta$ be a lattice homomorphism is strictly stronger than the condition that $\eta$ be order-preserving.

The fibers of a lattice homomorphism from~$L$ to another lattice constitute an equivalence relation $\equiv$ on $L$. An equivalence relation which arises in this way is called a {\em lattice congruence} on~$L$.
More directly, an equivalence relation $\equiv$ on~$L$ is a lattice congruence if and only if $a_1\equiv a_2$ and $b_1\equiv b_2$ implies $(a_1\join b_1)\equiv(a_2\join b_2)$ and $(a_1\meet b_1)\equiv(a_2\meet b_2)$.
It is an easy exercise to show that an equivalence relation $\Theta$ on a finite lattice~$L$ is a lattice congruence if and only if it satisfies the following three conditions (where $[x]_\Theta$ denotes the $\Theta$-equivalence class of~$x$):
\begin{enumerate}
\item[(i) ]Each equivalence class $[x]_\Theta$ is an interval in~$L$;
\item[(ii) ]The map $\pidown^{\Theta}$ taking~$x$ to the minimal element of $[x]_\Theta$ is order-preserving;\\[-2.5 mm]
\item[(iii) ]The map $\piup_{\Theta}$ taking~$x$ to the maximal element of $[x]_\Theta$ is order-preserving. 
\end{enumerate}

Given a lattice congruence $\Theta$ on~$L$, the {\em quotient lattice} $L/\Theta$ is the lattice whose elements are the congruence classes, with join and meet defined by $[x]_\Theta\join[y]_\Theta=[x\join y]_\Theta$ and $[x]_\Theta\meet[y]_\Theta=[x\meet y]_\Theta$.
Equivalently, $L/\Theta$ is the partial order on congruence classes which sets $[x]_\Theta\le[y]_\Theta$ if and only if there exists $x'\in[x]_\Theta$ and $y'\in[y]_\Theta$ such that $x'\le y'$ in~$L$.
It is easy to check that $L/\Theta$ is isomorphic to the subposet of~$L$ induced by the set $\pidown^{\Theta}(L)=\set{x\in L:\pidown^{\Theta}(x)=x}$.
Note that $\pidown^{\Theta}(L)$ need not be a sublattice of~$L$.
The following is \cite[Proposition~2.2]{con_app}.
\begin{prop}
\label{quotient covers}
Let~$L$ be a finite lattice, $\Theta$ a congruence on~$L$ and $x\in L$.
Then the map $y\mapsto[y]_\Theta$ restricts to a one-to-one correspondence between elements of~$L$ covered by $\pidown^\Theta(x)$ and elements of $L/\Theta$ covered by $[x]_\Theta$.
\end{prop}

We now remind the reader of some basic facts about Coxeter groups.
We also establish notation for what follows.
Throughout the paper,~$W$ denotes a finite Coxeter group of rank~$n$ with simple generators~$S$. For $J\subseteq S$, let $W_J$ be the subgroup of~$W$ generated by $J$, called a {\em standard parabolic subgroup} of~$W\!$.
Most often $J$ will be $S \setminus \set{s}$ for some \mbox{$s\in S$}; we write $\br{s}$ for $S \setminus \set{s}$.
The {\em reflections} of~$W$ are those elements which are conjugate to elements of~$S$.
The set of reflections is written~$T$.
An element $w\in W$ can be written as a word in~$S$.
A word for~$w$ which is minimal in length among words for~$w$ is called {\em reduced} and the {\em length} of~$w$, written $\ell(w)$, is the length of a reduced word for~$w$.
An {\em inversion} of~$w$ is a reflection $t\in T$ such that $\ell(tw)<\ell(w)$. If $s_1 s_2 \cdots s_{\ell}$ is a reduced word for $w$ then the inversions of $w$ are $s_1$, $s_1 s_2 s_1$, \dots, $s_1 s_2 \cdots s_{\ell} \cdots s_2 s_1$. 
The set of inversions is written $I(w)$ and~$w$ is uniquely determined by $I(w)$.
The {\em (right) weak order} on~$W$ is the partial order on~$W$ induced by containment of inversion sets.
Equivalently, the weak order on~$W$ is the transitive closure of the cover relations $w\covered ws$ whenever \mbox{$s\in S$} and $\ell(w)<\ell(ws)$.
This is further equivalent to the partial order defined by $v\le w$ if and only if there is a reduced word $a$ for~$w$ such that  some prefix (initial subword) of $a$ is a word for~$v$.
The weak order is known to be a lattice when~$W$ is finite.

The symbol~$W$ will now denote both the group~$W$ and the set~$W$ viewed as a poset (lattice).
All references to a partial order on~$W$ will refer to the weak order.
The phrase ``join-irreducible elements of the weak order on~$W$'' will be abbreviated to ``join-irreducibles of~$W$'' and similarly we will refer to ``meet-irreducibles of $W.$''

The unique maximal element of~$W$ is called $w_0$.
Conjugation by $w_0$ is an automorphism of the weak order and in particular permutes the simple generators~$S$.
The map $w\mapsto ww_0$ is an antiautomorphism of the weak order on $W,$ and in particular $I(ww_0)=T \setminus I(w)$.

A simple reflection~$s$ is called a \emph{descent} of~$w$ if $\ell(ws) < \ell(w)$ and an \emph{ascent} of~$w$ if $\ell(ws) > \ell(w)$. A {\em cover reflection} of $w\in W$ is a reflection~$t$ such that $tw\covered w$; the set of cover reflections of~$w$ can also be described as those reflections of the form $wsw^{-1}$ for~$s$ a descent of~$w$. The set of cover reflections is denoted by $\cov(w)$. There is one cover reflection of~$w$ for each element of~$W$ covered by~$w$.
Thus the join-irreducibles of~$W$ are the elements with exactly one cover reflection, or equivalently one descent. 
The following is \cite[Lemma~2.8]{sort_camb}.

\begin{lemma}
\label{cov w br s}
For $x\in W_{\br{s}}$, $\cov(s\join x)=\cov(x)\cup\set{s}$.
\end{lemma}

The map $x\mapsto sx$ is an involutive isomorphism between the intervals $[s,w_0]$ and $[1,sw_0]$.
In particular, we have the following observation, which we record now to avoid giving the simple argument repeatedly later.
\begin{lemma}
\label{s ji}
Let~$w$ be join-irreducible with $s<w$.
Then $sw$ is join-irreducible and if~$t$ is the unique cover reflection of~$w$ then $sts$ is the unique cover reflection of $sw$.
\end{lemma}
\begin{proof}
Since~$w$ is join-irreducible,~$w$ covers at most one element of $[s,w_0]$. If~$w$ does not cover any element of $[s,w_0]$ then $sw$ does not cover any element of $[1,s w_0]$ and must thus be~$1$. But this contradicts the assumption that $w \neq s$. Thus,~$w$ covers exactly one element of $[s,w_0]$, say $w \gtrdot wr$. Then $sw$ covers $swr$ and no other element of~$W\!$. The unique cover reflections of~$w$ and $sw$ are $wrw^{-1}$ and $sw r w^{-1} s$ respectively.
\end{proof}

For each $w\in W$ and each subset $J$ of~$S$ there is a unique factorization $w=w_J\cdot\closeleftq{J}{w}$ such that $w_J\in W_J$ and $\closeleftq{J}{w}$ satisfies $s\not\le\closeleftq{J}{w}$ for every $s\in J$.
The element $w_J$ appearing in this factorization is the unique element $w_J$ such that $I(w_J)=I(w)\cap W_J$.
For a fixed $w\in W,$ the set of elements~$x$ such that $x_J=w_J$ is an interval in $W,$ specifically the interval $[w_J,w_J\cdot\closeleftq{J}{(w_0)}]$. The map $w \mapsto w_J$ is a lattice homomorphism. 

For more on the factorization $w=w_J \cdot \closeleftq{J}{w}$, see Section~2.4 of~\cite{Bj-Br}. All the claims of the preceding paragraph except for the last one are either in~\cite{Bj-Br} or are easy consequences of results proved there.
The fact that $w \mapsto w_J$ is a lattice homomorphism is proven for example in~\cite{Jed} or \cite[Proposition~6.3]{congruence}.

\section{Sortable elements and Cambrian congruences}
\label{sort sec}
In this section we review definitions and quote or prove some preliminary results about sortable elements and Cambrian congruences.
For more details, see~\cite{sortable,sort_camb}.

A {\em Coxeter element} of~$W$ is an element of~$W$ of the  form $s_1s_2\cdots s_n$, where $s_1,s_2,\ldots, s_n$ are the simple generators $S,$ listed in any order.
(Recall that $n=|S|$.)
Two orderings of the generators produce the same Coxeter element if and only if they are related by a sequence of transpositions of adjacent generators which commute in~$W\!$.
A generator \mbox{$s\in S$} is {\em initial} in~$c$, or is an {\em initial letter of $c$}, if~$c$ can be written $s_1s_2\cdots s_n$ with $s_1=s$.
{\em Final letters} are defined similarly.

Given $w \in W,$ the half-infinite word
\[c^\infty=s_1s_2\cdots s_ns_1s_2\cdots s_ns_1s_2\cdots s_n\ldots\]
contains infinitely many subwords which are reduced words for~$w$.
The {\em $c$-sorting word} for $w\in W$ is the lexicographically leftmost subword of $c^\infty$ which is a reduced word for~$w$.
Inserting {\em dividers} ``$|$'' into $c^\infty$
\[c^\infty=s_1s_2\cdots s_n|s_1s_2\cdots s_n|s_1s_2\cdots s_n|\ldots,\]
we view the $c$-sorting word for~$w$ as a sequence of subsets of $S$, namely the sets of letters of the $c$-sorting word which occur between adjacent dividers.

An element $w\in W$ is {\em $c$-sortable} if its $c$-sorting word defines a sequence of subsets which is weakly decreasing under inclusion.
Formally, this definition requires a choice of reduced word for~$c$.
However, for a given~$w$, the $c$-sorting words for~$w$ arising from different reduced words for~$c$ are related by commutations of letters, with no commutations across dividers.
Thus in particular, the set of $c$-sortable elements does not depend on the choice of reduced word for~$c$.

\begin{example}
\label{B2sort}
Consider the Coxeter group $B_2$ with simple generators $s_0$ and~$s_1$.
When~$W$ is $B_2$ and $c=s_0s_1$, the $c$-sortable elements of~$W$ are~$1$,~$s_0$, $s_0s_1$, $s_0s_1s_0$, $s_0s_1s_0s_1$ and~$s_1$.  
The non-$c$-sortable elements are $s_1s_0$ and $s_1s_0s_1$.
\end{example}

The definition of sortability in terms of $c^{\infty}$ is intuitive but is not always the most helpful definition. The following two lemmas, which are \cite[Lemmas~2.4 and~2.5]{sortable}, give a recursive description of $c$-sortability.
\begin{lemma}
\label{sc}
Let~$s$ be an initial letter of~$c$ and let $w\in W$ with $s\not\le w$.
Then~$w$ is $c$-sortable if and only if it is an $sc$-sortable element of~$W_{\br{s}}$.
\end{lemma}

\begin{lemma}
\label{scs}
Let~$s$ be an initial letter of~$c$ and let $w\in W$ with $s\le w$.
Then~$w$ is $c$-sortable if and only if $sw$ is $scs$-sortable.
\end{lemma}

%To decide whether $s\le w$ means to decide whether the maximal cone of $\F$ corresponding to~$w$ is above the reflecting hyperplane for the reflection~$s$.
In~\cite{sortable}, the two lemmas above appear with the hypothesis $\ell(sw)>\ell(w)$ (resp.\ $\ell(sw)<\ell(w)$) instead of $s\not\le w$ (resp.\ $s\le w$).  
The characterization of the weak order in terms of inversion sets reconciles these two ways of stating the hypothesis.
In Lemma~\ref{sc}, $W_{\br{s}}$ is a Coxeter group of rank $n-1$ and in Lemma~\ref{scs}, $\ell(sw)<\ell(w)$, so these two lemmas characterize the $c$-sortable elements by induction on rank and length.
(The identity element~$1$ is $c$-sortable for any~$c$.)

For each~$c$, we define a map $\pidown^c$ from~$W$ to the $c$-sortable elements of~$W\!$.
The notation $\pidown^c$ suggests the order-theoretic characterization of lattice congruences given in Section~\ref{weak sec}.
For any Coxeter element~$c$, let $\pidown^c(1)=1$ and for~$s$ an initial letter of~$c$, define
\[
\pidown^c(w)=\left\lbrace\begin{array}{ll}
s\cdot\pidown^{scs}(sw)&\mbox{if }s\le w,\mbox{ or}\\
\pidown^{sc}(w_{\br{s}})&\mbox{if }s\not\le w.
\end{array}\right.
\]
In \cite[Section~3]{sort_camb}, it is shown that $\pidown^c(w)$ is the unique maximal $c$-sortable element weakly below~$w$.
Furthermore, it is shown that the fibers of $\pidown^c$ are a lattice congruence on~$W\!$, denoted by $\Theta_c$.
In particular, $\pidown^c$ is order preserving. 
In \cite[Section~5]{sort_camb}, $\Theta_c$ is identified as the {\em $c$-Cambrian congruence} on~$W$ in the sense of~\cite{cambrian}.
(Although we leave the lattice-theoretic details to~\cite{cambrian} and~\cite{sort_camb}, we will adopt the name ``$c$-Cambrian congruence'' for $\Theta_c$.)
We use the abbreviation $[w]_c$ for $[w]_{\Theta_c}$.

The {\em $c$-Cambrian lattice} is the quotient $W/\Theta_c$ of~$W$ modulo the $c$-Cambrian congruence.
This quotient is isomorphic to the restriction of the weak order to the $c$-sortable elements of~$W\!$.
Despite this isomorphism, to avoid confusion the notation of this paper will describe the $c$-Cambrian lattice as a partial order on congruence classes $[w]_c$, while comparisons of elements~$w$ will refer to the weak order on~$W\!$.
Thus $[w]_c\covered [w']_c$ is a cover relation in the quotient $W/\Theta_c$, while $w\covered w'$ says that~$w$ is covered by $w'$ in the weak order on~$W\!$.

\begin{figure}
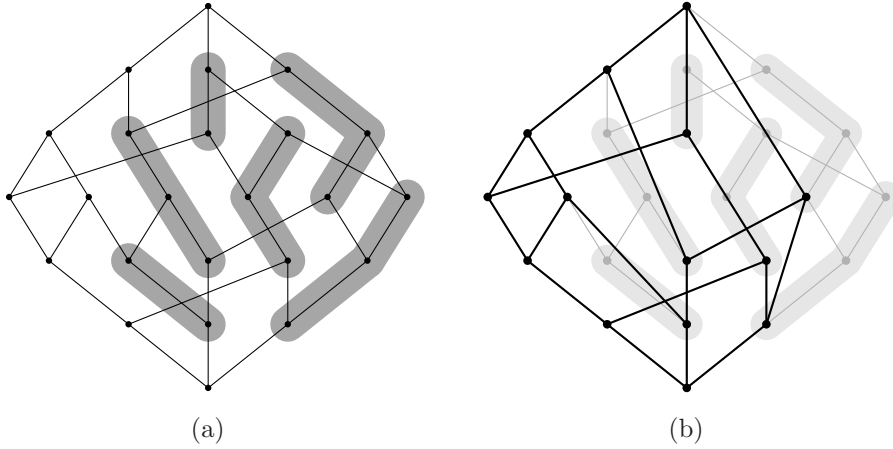

\centerline{
\begin{tabular}{ccc}
\scalebox{.4}{\epsfbox{weakS4cong.ps}}
&&
\scalebox{.4}{\epsfbox{weakS4bottoms.ps}}
\\[4 pt]
(a)&&(b)
\end{tabular}
}

\caption{A Cambrian congruence and the associated Cambrian lattice}
\label{A3cambfig}
\end{figure}

\begin{example}
\label{A3camb}
Figure~\ref{A3cambfig}.a shows the $s_1s_2s_3$-Cambrian congruence on the weak order for~$W$ of type $A_3$.
The gray shading indicates congruence classes of cardinality greater than one, and each unshaded vertex is a singleton congruence class.
The $s_1s_2s_3$-Cambrian lattice is the partial order on the congruence classes, as explained in Section~\ref{weak sec}.
Equivalently, the $s_1s_2s_3$-Cambrian lattice is the restriction of the weak order to $s_1s_2s_3$-sortable elements (bottom elements of congruence classes), as indicated in Figure~\ref{A3cambfig}.b.
\end{example}

The $c$-Cambrian congruence has an upward projection map $\piup_c$ which takes each $w\in W$ to the top element of its $\Theta_c$-congruence class.
This map is given by $\piup_c(w)=\left(\pidown^{(c^{-1})}(ww_0) \right)w_0$, and satisfies the following recursion when~$s$ is \emph{final} in $c$: 
\[
\piup_c(w)=\left\lbrace\begin{array}{ll}
s\cdot\piup_{scs}(sw)&\mbox{if }s\not\le w,\mbox{ or}\\
\piup_{cs}(w_{\br{s}})\cdot\leftq{\br{s}}{w_0}&\mbox{if }s\le w.
\end{array}\right.
\]
In particular the antiautomorphism $w\mapsto ww_0$ takes $c$-Cambrian congruence classes to $(c^{-1})$-Cambrian congruence classes.
An element $w\in W$ is called {\em $c$-antisortable} if $\piup_c(w)=w$.
Equivalently,~$w$ is $c$-antisortable if and only if $ww_0$ is $(c^{-1})$-sortable.
The map $\piup_c$ takes~$w$ to the unique minimal $c$-antisortable element above~$w$.

The elements~$1$ and $w_0$ are alone in their $c$-Cambrian congruence classes.
An initial letter~$s$ of~$c$ also constitutes a singleton congruence class.

We now record three simple lemmas about $\pidown^c$, $\piup_c$ and $c$-antisortable elements.

\begin{lemma}
\label{no crossing}
Let~$c$ be a Coxeter element, let~$r$ be a simple reflection and let $w\in W$. Then $r\le w$ if and only if $r \le \pidown^c(w)$. 
\end{lemma}

\begin{proof}
First, suppose that $\pidown^c(w) \geq r$. 
Since $\pidown^c(w) \leq w$, we also have $w \geq r$. 
Now, suppose $w \geq r$. Then, since $\pidown^c$ is order preserving, $\pidown^c(w) \geq \pidown^c(r)=r$.
\end{proof}

\begin{lemma}
\label{scs dual}
Let~$s$ be an initial letter of~$c$ and let $w\in W$ with $s\le w$.
Then~$w$ is \mbox{$c$-antisortable} if and only if $sw$ is $scs$-antisortable.
\end{lemma}
\begin{proof}
The element~$w$ is $c$-antisortable if and only if $ww_0$ is $(c^{-1})$-sortable.
Observe that~$s$ is initial in $sc^{-1}s$ and that $s\le sww_0$. Thus, by Lemma~\ref{scs},  $sww_0$ is $sc^{-1}s$-sortable if and only if $ww_0$ is $(c^{-1})$-sortable.
Now $sww_0$ is $sc^{-1}s$-sortable if and only if $sw$ is $scs$-antisortable.
\end{proof}

\begin{lemma}
\label{piup s}
If~$s$ is final in~$c$ then $\piup_c(s)=w_0\cdot\left((w_0)_{\br{s'}}\right)$ for $s'=w_0sw_0$.
\end{lemma}
\begin{proof}
The reflection $s'$ is a simple reflection because $w_0$ permutes~$S$.
By the recursive characterization, $\piup_c(s)$ is equal to
\[\piup_{cs}(1)\cdot\leftq{\br{s}}{w_0}=\leftq{\br{s}}{w_0}=(w_0)_{\br{s}}\cdot w_0=
w_0\cdot\left(w_0(w_0)_{\br{s}}w_0\right)=w_0\cdot\left((w_0)_{\br{s'}}\right). \qedhere \]
\end{proof}

The proofs in this paper will also rely on nontrivial properties of sortable elements which we now quote.
We begin with \cite[Theorem~6.1]{sortable}.
In the present paper, it will not be necessary to define noncrossing partitions or the map $\nc_c$.
Full details, including citations, are found in \cite[Sections~5 and~6]{sortable}.
Following the statement of the theorem we discuss how it applies to the present context.
\begin{theorem}
\label{nc}
For any Coxeter element~$c$, the map $w\mapsto\nc_c(w)$ is a bijection from the set of $c$-sortable elements to 
the set of noncrossing partitions with respect to~$c$.
Furthermore, $\nc_c$ maps $c$-sortable elements with~$k$ descents to $c$-noncrossing partitions of rank~$k$.
\end{theorem}
The noncrossing partitions (with respect to~$W$ and~$c$) of rank~$1$ are exactly the reflections in~$W\!$.
For the purposes of this section, all we need to know about $\nc_c$ is Theorem~\ref{nc} and the following fact: 
If~$w$ is a $c$-sortable element with~$1$ descent (i.e.\ a $c$-sortable join-irreducible) then $\nc_c(w)$ is the unique cover reflection of~$w$. Thus we have the following corollary to Theorem~\ref{nc}.
\begin{cor}
\label{irr cov}
For each reflection~$t$ of $W,$ there is exactly one $c$-sortable join-irreducible whose unique cover reflection is~$t$.
\end{cor}

The number of noncrossing partitions of rank $n-1$ is also equal to the number~$|T|$ of reflections in~$W\!$.
Thus $c$-sortable meet-irreducibles are also counted by $|T|$. (Recall that an element of $W$ is meet-irreducible if and only if it is covered by exactly one element, or equivalently, covers exactly $n-1$ elements.)
Because the map $w\mapsto ww_0$ is an antiautomorphism of the weak order on~$W$ and takes $c$-sortable elements to $(c^{-1})$-antisortable elements, the same is true of $c$-antisortable join-irreducible or meet-irreducible elements.
We summarize in the following corollary to Theorem~\ref{nc}.

\begin{corollary}
\label{num irr}
For~$W$ a finite Coxeter group and~$c$ a Coxeter element of $W,$ the following numbers are all equal:
\begin{enumerate}
\item[(i) ] the number of $c$-sortable join-irreducibles of~$W$;
\item[(ii) ] the number of $c$-sortable meet-irreducibles of~$W$;
\item[(iii) ] the number of $c$-antisortable join-irreducibles of~$W$;
\item[(iv) ] the number of $c$-antisortable meet-irreducibles of~$W$; and 
\item[(v) ] the number of reflections in~$W\!$.
\end{enumerate}
\end{corollary}

We conclude the section by quoting \cite[Theorem~1.2]{sort_camb} and using it to prove a lemma.
\begin{theorem}
\label{sublattice}
Let~$c$ be a Coxeter element of a finite Coxeter group~$W\!$.\@
The \mbox{$c$-sortable} elements constitute a sublattice of the weak order on $W\!$.\@
\end{theorem}

\begin{lemma}
\label{anti cov}
Let~$s$ be an initial letter of~$c$.
If~$w$ is $c$-antisortable and $s\not\le w$ then 
\begin{enumerate}
\item[(i) ]$sw\covers w$, 
\item[(ii) ]$s\join w=sw$ and
\item[(iii) ]$w$ is $scs$-antisortable.
\end{enumerate}
\end{lemma}
\begin{proof}
Let~$w$ be~$c$-antisortable with $s\not\le w$.
The set of elements~$x$ such that $x_{\br{s}}=w_{\br{s}}$ is an interval $I$ in~$W$. Now, if $x \in I$ and $x \not \geq s$ then $\pidown^c(x)=\pidown^{sc}(x_{\br{s}})=\pidown^{sc}(w_{\br{s}})=\pidown^c(w)$ so, by the assumption that $w$ is $c$-antisortable, we have $w \geq x$ for such an $x$.

On the other hand, since $(s\join w)_{\br{s}} = s_{\br{s}} \join w_{\br{s}}=w_{\br{s}}$, $I$ contains elements which are greater than or equal to~$s$. 
In particular,~$w$ is not maximal in $I$ so we take $w'$ to be some element covering~$w$ with $(w')_{\br{s}}=w_{\br{s}}$. By the observation in the first paragraph, $w' \geq s$. But since $w' \covers w$ with $w \not \geq s$ and $w' \geq s$, we must have $ sw=w'\covers w$.
Furthermore, $sw=s\join w$ because $s\le sw$. By the dual of Theorem~\ref{sublattice}, $sw$ is $c$-antisortable and, by Lemma~\ref{scs dual},~$w$ is $scs$-antisortable.
\end{proof}

\section{Cambrian fans}
\label{camb fan sec}
In this section we define the $c$-Cambrian fan for each finite Coxeter group~$W$ and Coxeter element~$c$ of~$W\!$.
We also prove a few preliminary results.
In Section~\ref{ray sec} we make a careful study of the rays of the $c$-Cambrian fan.

An {\em arrangement} of hyperplanes in a vector space~$V$ is a collection of hyperplanes (codimension~$1$ subspaces). A central arrangement is called {\em central} if all of the hyperplanes pass through the origin.  
That is, the hyperplanes are linear subspaces rather than affine subspaces.
A central arrangement is called {\em essential} if the intersection of the hyperplanes is the origin.

We continue to let $(W,S)$ be a finite Coxeter system of rank~$n$ with reflections~$T$ and longest element $w_0$. 
Fix some root system $\Phi$ for~$W$ and let~$V$ be the geometric representation of~$W$; we write $V(W)$ when it is necessary for clarity.
This is the representation of~$W$ on the real vector space spanned by the root system of~$W\!$.
A reflection~$t$ of~$W$ acts by the orthogonal reflection
\[v \mapsto v - 2\frac{\br{\alpha_t, v}}{\br{\alpha_t,\alpha_t}} \alpha_t,\]
where $\alpha_t$ is the positive root corresponding to~$t$ and $\br{\cdot,\cdot}$ is the usual inner product.
The hyperplane fixed by the reflection~$t$ is denoted by~$H_t$. 
The {\em Coxeter arrangement}~$\A$ for~$W$ is the collection of all such hyperplanes; we will write $\A(W)$ when necessary. The complement $V \setminus (\bigcup\A)$ of~$\A$ is composed of open cones whose closures are called {\em regions}.
The regions are in canonical bijective correspondence with the elements of $W,$ and each region has~$n$ facet hyperplanes. More specifically, the \emph{dominant chamber} $D:=\bigcap_{s \in S} \{ v : \langle v,\alpha_s \rangle \geq 0 \}$ corresponds to the identity and $wD$ corresponds to~$w$. 

A subset $U$ of~$V$ is {\em below} a hyperplane $H\in\A$ if every point in $U$ is either on~$H$ or is on the same side of~$H$ as $D$.
The subset is {\em strictly below}~$H$ if it is below~$H$ and does not intersect~$H$.
Similarly, $U$ can be {\em above} or {\em strictly above}~$H$.
The inversions of an element $w\in W$, defined in Section~\ref{weak sec} to be those reflections $t$ for which $\ell(tw) <   \ell(w)$, can also be described as the reflections $t$ such that $wD$ is above $H_t$. This result is perhaps more frequently quoted in its dual form: a reflection $t$ is an inversion of $w$ if and only if $w^{-1}(\alpha_t)$ is a negative root \cite[Proposition~4.4.6]{Bj-Br}.
In the case of a simple reflection \mbox{$s\in S$}, $\ell(sw)<\ell(w)$ if and only if $s\le w$ in the weak order. 
Thus deciding whether $wD$ is above $H_s$ or below $H_s$ is a weak order comparison.

The Coxeter arrangement~$\A$ is a central, essential arrangement. 
If $J \subset S$, let $\A_J$ be the subset of~$\A$ consisting of those hyperplanes $H_t \in \A(W)$ for which $t$ is in $W_J$. 
Then $\A_J$ is a central arrangement but it is not essential. 
If $I_J=\bigcap_{H \in \A_J} H$  then we have $V(W)\cong I_J \times V(W_J)$ and each hyperplane~$H_t$ in $\A_J$ is the direct product of $I_J$ with the hyperplane~$H_t$ in $\A(W_J)$.
We write $\Proj_J$ for the linear projection $V(W) \onto V(W) / I_J \isomorph V(W_J)$. 
This projection will be used in Lemma~\ref{para ray} and the proof of Theorem~\ref{cl iso}. The following proposition relates the geometric projection $\Proj_J$ to the combinatorial projection $w \mapsto w_J$.
\begin{prop} \label{projection}
For $w \in W$, we have $\Proj_J(w D) \subseteq w_J D_J$, where $D_J$ is the dominant chamber for $\A(W_J)$.
\end{prop}

A {\em fan} $\F$ is a family of nonempty closed polyhedral (and in particular convex) cones in~$V$ such that
\begin{enumerate}
\item[(i) ] For any cone in $\F,$ all faces of that cone are also in $\F,$ and
\item[(ii) ] The intersection of two cones in $\F$ is a face of both.
\end{enumerate}
A fan is {\em complete} if its union is all of~$V\!$.
It is {\em essential} (or {\em pointed}) if the intersection of all of the cones of $\F$ is the origin.
For more information about fans, see \cite[Lecture~7]{Ziegler}.

Let $\F$ be the fan consisting of the regions of~$\A$ and all of their faces.
The fan $\F$ is complete and essential. (See \cite[Sections~1.12--1.15]{Humphreys}.)
The faces of $\F$ have an elegant description: they are in bijection with pairs $(w,J)$ where~$w$ is an element of~$W$ and $J$ is a subset of the ascents of~$w$. 
The pair $(w,J)$ corresponds to $C(w,J):=w \cdot \left( D \cap \bigcap_{s \in J} \{ v : \langle v,\alpha_s \rangle = 0 \} \right)$. 
We may recover~$w$ from $C(w,J)$ by the fact that~$w$ is the smallest element of~$W$ (in weak order) such that $wD$ contains $C(w,J)$.
It is then easy to recover $J$. The cone $C(w,J)$ has dimension $n-|J|$.

Let $\G$ be a complete fan in $\RR^n$ and let~$v$ be a generic vector in $\RR^n$. 
Suppose that the intersection of two maximal cones~$C$ and~$C'$ spans a hyperplane~$H$ with~$v$ on the same side of~$H$ as~$C$.
We put $C' > C$. In general, it is possible that there is a sequence $C_1$, $C_2$, \ldots, $C_r$ of maximal cones of $\G$ such that $C_1 > C_2 > \cdots > C_r > C_1$. If this does not occur, then we define a poset on the maximal cones of $\G$ by taking the transitive closure of all relations $C' > C$ and we say that this poset is {\em induced\,}\footnote{An analogous construction in~\cite{con_app} featured posets induced on fans by linear functionals.  The vector~$v$ occurring here points in the direction which minimizes the linear functional of~\cite{con_app}.}
 on $\G$ by~$v$; roughly speaking, going ``down'' in the poset means moving in the direction of~$v$. So, for example, the weak order is induced on $\F$ by any $v$ in the interior of $D$.  
If~$C$ is a maximal cone of $\G$ and if $C$ is simplicial, then we define the {\em bottom face} of~$C$ with respect to~$v$ to be the minimal (under containment) face $F$ of~$C$ such that for any vector~$x$ in the relative interior of $F$, there exists an $\ep>0$ such that $x-\ep v$ is in~$C$. 
In other words, the bottom face of~$C$ is the intersection of the facets $F$ of~$C$ which separate~$C$ from a face lower than~$C$ in the poset induced by~$v$.

%If~$C$ is a simplicial cone of $\G$ then we define the {\em bottom face} of~$C$ with respect to~$v$ to be the minimal (under containment) face $F$ of~$C$ such that for any vector~$x$ in the relative interior of $F$, there exists an $\ep>0$ such that $x-\ep v$ is in~$C$. 
%If~$C$ is a maximal cone then the bottom face of~$C$ is the intersection of the facets $F$ of~$C$ which separate~$C$ from a face lower than~$C$ in the poset induced by~$v$.

We now review a construction from~\cite{con_app} which, given an arbitrary lattice congruence $\Theta$ on the weak order on $W,$ constructs a complete fan~$\F_\Theta$ which coarsens $\F$ (in the sense that every cone of~$\F_\Theta$ is a union of cones of $\F$).
The maximal cones of~$\F_\Theta$ correspond to congruence classes of $\Theta$.
Specifically, each maximal cone is the union of the regions of~$\A$ corresponding to the elements of the congruence class.
In~\cite[Section~5]{con_app} it is shown that the collection~$\F_\Theta$ consisting of these maximal cones together with all of their faces is indeed a complete fan.  
In what follows, we identify a congruence class with the corresponding maximal cone of~$\F_\Theta$.

\begin{example}
\label{B2FcExample}
For $W=B_2$ with $c=s_0s_1$, the fan $\F$ is shown in Figure~\ref{B2Fcfig}.a, with maximal cones labeled by elements of~$W\!$.
Figure~\ref{B2Fcfig}.b shows~$\F_c$, with maximal cones labeled by $c$-sortable elements.
(Cf.\ Example~\ref{B2sort}.)
The weak order on $B_2$ is the poset on the regions of $\F$ such that one moves up in the partial order by passing to an adjacent region which is ``higher'' on the page.
The $c$-Cambrian lattice is the poset on the maximal cones of~$\F_c$ with a similar description.
\end{example}

\begin{figure}[ht]
\centerline{
\begin{tabular}{cc}
\scalebox{1}{
\epsfbox{B2F.ps}
\begin{picture}(0,0)(75,-75)
\put(-6,-60){$1$}
\put(-51,-45){$s_0$}
\put(-66,-2){$s_0s_1$}
\put(-56,35){$s_0s_1s_0$}
\put(-20,60){$s_0s_1s_0s_1$}
\put(35,-45){$s_1$}
\put(40,-2){$s_1s_0$}
\put(25,35){$s_1s_0s_1$}
\end{picture}
}
&
\scalebox{1}{
\epsfbox{B2Fc.ps}
\begin{picture}(0,0)(75,-75)
\put(-6,-60){$1$}
\put(-51,-45){$s_0$}
\put(-66,-2){$s_0s_1$}
\put(-56,35){$s_0s_1s_0$}
\put(-20,60){$s_0s_1s_0s_1$}
\put(24,-2){$s_1$}
\end{picture}
}
\\[4 pt]
(a)&(b)
\end{tabular}
}
\caption{The fans $\F$ and~$\F_c$}
\label{B2Fcfig}
\end{figure}

%
%\begin{example}
%\label{A3abc}
%Figure~\ref{A3abcfig} shows the fan~$\F_\Theta$ for the congruence $\Theta$ illustrated in Figure~\ref{A3cambfig}.a.
%(This is $\Theta_c$ for $c=s_1s_2s_3$ and~$W$ of type $A_3$.)
%The solid lines indicate the intersection of~$\F_\Theta$ with a unit sphere.
%The dotted lines indicate how each maximal cone of~$\F_\Theta$ is partitioned into maximal cones of the Coxeter fan.
%Black lines lie on the hemisphere closer to the viewer and grey lines lie on the far hemisphere.
%The maximal cone for~$1$ (the minimal element of the weak order) is at the center of the near hemisphere.
%\end{example}

%
%\begin{figure}
%\centerline{\scalebox{.35}{\epsfbox{A3abc.ps}}}
%\caption{A $c$-Cambrian fan}
%\label{A3abcfig}
%\end{figure}

The lattice $W/\Theta$ is a partial order on the maximal cones of~$\F_\Theta$.
In fact, the pair ($\F_\Theta,W/\Theta)$ is a {\em fan poset} \cite[Theorem~1.1]{con_app}, and thus by \cite[Proposition~3.3]{con_app}, we have the following.

\begin{prop} \label{ItsACover}
Let $[w]_\Theta$ and $[w']_\Theta$ be maximal cones of~$\F_\Theta$. 
Then $[w]_\Theta$ and $[w']_\Theta$ are a covering pair in $W/\Theta$ if and only if they intersect in a common facet.
\end{prop}

The proof of the following lemma is essentially contained in the proof of \cite[Proposition~5.5]{con_app}.
However, since that result is stated quite differently and in broader generality, we give a proof here.
The dual statement about the upward projection $\piup_{\Theta}$ also holds.

\begin{lemma}
\label{separators}
Let $\pidown^{\Theta}$ be the downward projection map associated to a lattice congruence on the weak order on a finite Coxeter group.
A hyperplane~$H$ separates a congruence class $[w]_\Theta$ from a congruence class $[x]_\Theta\covered [w]_\Theta$ if and only if~$H$ separates $\pidown^{\Theta}(w)$ from an element of~$W$ covered by $\pidown^{\Theta}(w)$.
\end{lemma}
\begin{proof}
Let $\EL$ be the set of hyperplanes~$H$ in~$\A$ such that~$H$ separates the congruence class $[w]_\Theta$ from a congruence class covered by $[w]_\Theta$.
Since the congruence classes are convex cones, no two hyperplanes in $\EL$ separate $[w]_\Theta$ from the same congruence class covered by $[w]_\Theta$.
By Proposition~\ref{quotient covers}, congruence classes covered by $[w]_\Theta$ are in one-to-one correspondence with elements~$x$ covered by $\pidown^{\Theta}(w)$.
Each such~$x$ is separated from $\pidown^{\Theta}(w)$ by a distinct hyperplane in $\EL$, so $\EL$ is the set of hyperplanes~$H$ such that~$H$ separates $\pidown^{\Theta}(w)$ from an element covered by $\pidown^{\Theta}(w)$.
\end{proof}

The {\em $c$-Cambrian fan}~$\F_c$ is the essential fan $\F_{\Theta_c}$ arising from this construction, where $\Theta_c$ is the $c$-Cambrian congruence described in Section~\ref{sort sec}.
The $c$-Cambrian fan~$\F_c$ and the $c$-Cambrian lattice $W/\Theta_c$ have many pleasant properties following from a general theorem \cite[Theorem 1.1]{con_app} which applies to fans~$\F_\Theta$ and quotients $W/\Theta$ for general lattice congruences $\Theta$ on~$W\!$.
We list some of these properties here for emphasis.

\begin{enumerate}
\item[(i) ] Any linear extension of the $c$-Cambrian lattice is a shelling order of~$\F_c$.
\item[(ii) ] The $c$-Cambrian lattice is the order induced on the maximal cones of~$\F_c$ by any vector lying in the interior of $D$. 
\item[(iii) ] For any interval in the $c$-Cambrian lattice, the union of the corresponding cones of~$\F_c$ is a convex cone.
\item[(iv) ] For any cone $F$ in~$\F_c$, the set of maximal cones in~$\F_c$ containing $F$ is an interval in the $c$-Cambrian lattice.
\item[(v) ] A closed interval $I$ in the $c$-Cambrian lattice has proper part homotopy equivalent to an $(n-k-2)$-dimensional sphere if and only if there is some $k$-dimensional cone $F$ of~$\F_c$ such that $I$ is the set of all maximal 
cones of~$\F_c$ containing $F$.
\item[(vi) ] A closed interval $I$ has proper part homotopy equivalent to a $(k-2)$-dimensional sphere if and only if $I$ has~$k$ atoms and the join of the atoms of $I$ is the top element of $I$.
\item[(vii) ] If the proper part of a closed interval $I$ is not homotopy-spherical then it is contractible.
\end{enumerate}

By applying Lemma~\ref{separators} to the case of the $c$-Cambrian fan and appealing to Corollary~\ref{num irr}, we obtain the following useful fact about $c$-sortable and $c$-antisortable join-irreducibles.
Dually, the analogous statement for meet-irreducibles also holds.
\begin{prop}
\label{irr bij}
The upward projection $\piup_c$ restricts to a bijection from $c$-sortable join-irreducibles to $c$-antisortable join-irreducibles.
The inverse is the restriction of~$\pidown^c$.
\end{prop}
\begin{proof}
Let~$v$ be a $c$-sortable join-irreducible.
Then $v=\pidown^c(v)$ and since~$v$ covers exactly one element, by Lemma~\ref{separators} there is exactly one hyperplane separating $[v]_c$ from a congruence class covered by $[v]_c$.
Since $[v]_c$ is a full-dimensional cone in an essential fan, it must have at least~$n$ facet hyperplanes.
Thus (using Proposition~\ref{separators}) there must be at least $n-1$ congruence classes covering $[v]_c$ and by the dual of Lemma~\ref{separators} this means that $\piup_c(v)$ is covered by at least $n-1$ distinct elements.
The only element of~$W$ covered by~$n$ elements is~$1$, and every other element of~$W$ is covered by fewer elements.
If $\piup_c(v)$ is~$1$ then~$v$ is~$1$, contradicting the assumption that~$v$ is join-irreducible.
Therefore $\piup_c(v)$ is covered by exactly $n-1$ elements and covers exactly one element.
We have shown the $\piup_c$ maps $c$-sortable join-irreducibles to $c$-antisortable join-irreducibles.
Since $v=\pidown^c(\piup_c(v))$, the map is one-to-one and thus by Corollary~\ref{num irr} it is a bijection with inverse $\pidown^c$.
\end{proof}

We now describe the faces of~$\F_\Theta$.  (Cf.\ \cite[Proposition 5.10]{con_app}.)
We will use this description in Section~\ref{ray sec} when we discuss the rays of $\F_{c}$.
Let $w$ be maximal in its $\Theta$-equivalence class, in other words, let $w=\piup_{\Theta}(w)$, and let $J$ be a collection of ascents of~$w$.
In particular, for all $s \in J$, $[ws]_{\Theta} \neq [w]_{\Theta}$.
Recall that $C(w,J)$ is an $(n-|J|)$-dimensional face of the Coxeter fan $\F,$ and define $C_\Theta(w,J)$ to be the unique $(n-|J|)$-dimensional face of~$\F_\Theta$ containing $C(w,J)$.

To see that $C_\Theta(w,J)$ is well-defined, notice that by the dual of Lemma~\ref{separators}, there is a collection of $|J|$ facets of $[w]_\Theta$ all of which contain $C(w,J)$.
Furthermore, since each of these facets of $[w]_\Theta$ contains a facet of the region for~$w$ (a simplicial cone), the intersection of these facets is $(n-|J|)$-dimensional.
This is a face of~$\F_\Theta$ because it is an intersection of faces of~$\F_\Theta$.
Uniqueness is assured because no two distinct $k$-dimensional faces of a fan have a $k$-dimensional intersection.

\begin{prop} \label{faces of congruence}
The map $(w,J)\mapsto C_\Theta(w,J)$ is a bijection from ordered pairs $(w,J)$ with $w=\piup_{\Theta}(w)$ and $[ws]_\Theta\covers[w]_\Theta$ for every $s\in J$ to faces of~$\F_\Theta$.
\end{prop}
\begin{proof}
We only give full details of the proof for the restriction of the map to pairs such that $|J|=n-1$, i.e.\ $C_\Theta(w,J)$ is a ray.
Only that restriction is used in this paper.
For more general $J$, we sketch how a proof can be constructed using ideas, terminology and results of~\cite{con_app}.\footnote{In fact, this proposition, and the proof sketched here, is valid in the more general setting of \cite[Section~5]{con_app}.
We can replace the weak order~$W$ with a poset of regions of a simplicial hyperplane arrangement.
Instead of pairs $(w,J)$, we take pairs $(w,P)$ where~$w$ is a region maximal in $[w]_\Theta$ and $P$ is a set of facet hyperplanes of~$w$ separating $[w]_\Theta$ from classes above $[w]_\Theta$.}
The full proof is not difficult but would require quoting~\cite{con_app} in more detail than is desirable.

We first show that the map $(w,J) \to C_{\Theta}(w,J)$ is surjective. Let~$C$ be a $k$-dimensional face in~$\F_\Theta$.
By \cite[Theorem~1.1]{con_app} (cf.\ property (iv) of $c$-Cambrian fans, above), the set of maximal cones of~$\F_\Theta$ which contain~$C$ is an interval $I$ in $W/\Theta$.
Thus $I$ has a unique minimal element; this minimal element is a $\Theta$-congruence class $[w]_\Theta$, where~$w$ is chosen to be the maximal element in the class.
(In other words, $\piup_\Theta(w)=w$.)

Since $[w]_\Theta$ is lowest in $W/\Theta$ among congruence classes containing~$C$, every facet of $[w]_\Theta$ containing~$C$ separates $[w]_\Theta$ from a class that is higher in $W/\Theta$.
Thus, by Proposition~\ref{ItsACover}, every facet of $[w]_\Theta$ containing~$C$ separates $[w]_\Theta$ from a class $[ws]_\Theta\covers[w]_\Theta$. 
The hyperplanes separating~$w$ from those elements of~$W$ which cover it are transverse, and by Lemma~\ref{separators} these are exactly the hyperplanes separating $[w]_{\Theta}$ from the classes that cover it. 
Thus, since~$C$ is $k$-dimensional, there are exactly $n-k$ facets of $[w]_{\Theta}$ containing~$C$.
Let $J$ be the set of generators~$s$ such that $[ws]_\Theta\covers[w]_\Theta$ and such that the facet separating the two contains~$C$. 
Then $C(w,J)$ is $k$-dimensional and is contained in~$C$, so that $C=C_\Theta(w,J)$.
We have shown that the map is surjective.

The restriction of the map to pairs $(w,J)$ with $|J|=n-1$ is injective because a ray of~$\F_\Theta$ can't contain two distinct rays of $\F$.
We now sketch a proof that the unrestricted map is injective, continuing the notation of the previous two paragraphs.
Suppose $(w',J')$ obeys the conditions that $w'=\piup_{\Theta}(w')$, $[w's]_\Theta\covers[w']_\Theta$ for every $s\in J'$, $|J'|=|J|$ and $C(w',J')\subseteq C$.
We will show that $(w',J')=(w,J)$.
If $w'=w$ then since $J$ was defined by considering the set of all facets of $[w]_\Theta$ containing~$C$, we must have $J'=J$.

If $w'\neq w$ then $[w']_\Theta$ is in particular not minimal in $I$ with respect to $W/\Theta$.
Arguing as in the proof of \cite[Proposition~5.3]{con_app}, one shows that the interval $I$ is isomorphic to the quotient of a facial interval of $\F$ modulo the restriction of $\Theta$.
This restriction is bisimplicial by \cite[Proposition~5.5]{con_app}, leading to the conclusion that no element of $I$ (except the minimal element) is covered by $n-k$ or more elements of $I$.
Thus there do not exist $n-k$ distinct facets of $[w']_\Theta$ containing~$C$ and separating $[w']_\Theta$ from classes higher than $[w']_\Theta$ in $W/\Theta$.
In particular, the intersection of the facets separating $[w']_\Theta$ from $[w's]_\Theta$ for $s\in J'$ is a $k$-dimensional face of $[w']_\Theta$ distinct from~$C$.
This contradicts the supposition that $C(w',J')\subseteq C$, thus proving that $w'=w$ and thus that $(w',J')=(w,J)$.
\end{proof} 

\section{The cluster complex}
\label{cluster sec}
In this section we review the definition of the $c$-cluster complex, describe the map $\cl_c$ from the $c$-Cambrian fan to the $c$-cluster fan and give examples.
We begin by reviewing the definition of clusters in the sense of Fomin and Zelevinsky~\cite{ga} (as extended by Marsh, Reineke and Zelevinsky~\cite{MRZ} and extended slightly further in~\cite{sortable}).
Let~$\Phi$ be a root system for~$W$ with positive roots $\Phi_+$ and simple roots~$\Pi$.
For any reflection~$t$ of $W,$ let $\alpha_t$ denote the positive root associated to~$t$.
The roots in $\Pge=\Phi_+\cup(-\Pi)$ are called {\em almost positive roots}.
For any $J\subseteq S$, the set $(\Phi_J)_{\ge -1}$ is the intersection of $\Pge$ with the subset of~$\Phi$ corresponding to the parabolic subgroup $W_J$. 

For each \mbox{$s\in S$}, define an involution $\sigma_s:\Pge\to\Pge$ by
\[\sigma_s(\alpha):=\left\lbrace\begin{array}{ll}
\alpha&\mbox{if }\alpha\in(-\Pi)\mbox{ and }\alpha\neq-\alpha_s,\mbox{ or}\\
s(\alpha)&\mbox{otherwise.}
\end{array}\right.\]
The {\em $c$-compatibility} $\cm_c$ relation on $\Pge$ is defined by the following properties:
\begin{enumerate}
\item[(i) ]For any \mbox{$s\in S$}, $\beta\in\Pge$ and Coxeter element~$c$,
\[-\alpha_s\cm_c\beta\mbox{ if and only if }\beta\in(\Phi_{\br{s}})_{\ge-1}.\]
\item[(ii) ]For any $\alpha_1,\alpha_2\in\Pge$ and any initial letter~$s$ of~$c$, 
\[\alpha_1\cm_c \alpha_2\mbox{ if and only if }\sigma_s(\alpha_1)\,\cm_{scs}\,\sigma_s(\alpha_2).\]
\end{enumerate}
The relations $\cm_c$ and $\cm_{c^{-1}}$ coincide.
(See \cite[Proposition~3.1]{MRZ} and \cite[Proposition~7.4]{sortable}.)

A {\em $c$-compatible subset} of $\Pge$ is a set of pairwise $c$-compatible almost positive roots.
A {\em $c$-cluster} is a maximal $c$-compatible subset.
All $c$-clusters have cardinality~$n$.
Since each element of a $c$-cluster is a vector, the positive real span of the elements of a $c$-cluster is a well-defined cone.
In fact each $c$-cluster defines an $n$-dimensional cone, and these cones are the maximal cones in a complete fan (defined on the linear span of~$\Phi$).
This is the $c$-cluster fan.  
A set $\{\alpha_1,\ldots, \alpha_n \}$ is a $c$-cluster if and only if $\{ \sigma_s(\alpha_1), \ldots, \sigma_s(\alpha_n) \}$ is an $scs$-cluster. 
Thus, there is a continuous piecewise linear isomorphism between the $c$-cluster fan and the $scs$-cluster fan which is linear on each cone and sends each~$\alpha$ to $\sigma_s(\alpha)$.

\begin{figure}[ht]
\centerline{
\scalebox{1}{
\epsfbox{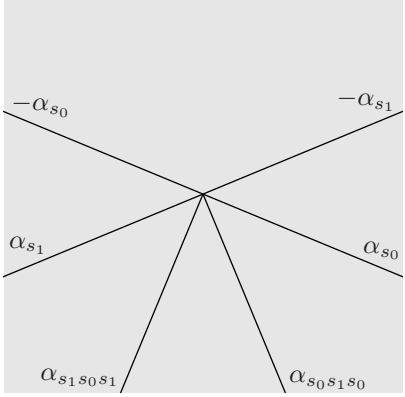}
\begin{picture}(0,0)(75,-75)
\put(-76,32){$-\alpha_{s_0}$}
\put(47,34){$-\alpha_{s_1}$}
\put(57,-22){$\alpha_{s_0}$}
\put(-77,-20){$\alpha_{s_1}$}
\put(29,-71){$\alpha_{s_0s_1s_0}$}
\put(-65,-70){$\alpha_{s_1s_0s_1}$}
\end{picture}
}}
\caption{The $c$-cluster fan}
\label{B2clusfan}
\end{figure}

\begin{example}
\label{B2clusfanEx}
The reflections in $W=B_2$ are $s_0$,~$s_1$, $s_0s_1s_0$ and $s_1s_0s_1$.
For $c=s_0s_1$, the $c$-cluster fan is shown in Figure~\ref{B2clusfan}, with each ray labeled by the corresponding almost positive root. 
\end{example}

The map $\cl_c$ takes a $c$-sortable element~$w$ to a set of~$n$ almost positive roots.
Let $a=a_1a_2\cdots a_k$ be the $c$-sorting word for~$w$.
If \mbox{$s\in S$} occurs in $a$ then the \emph{last reflection} for~$s$ in~$w$ is $a_1a_2\cdots a_{j-1} a_ja_{j-1}\cdots a_2a_1$, where $a_j$ is the rightmost occurrence of~$s$ in $a$.
The set $\cl_c(w)$ is obtained by taking the set of all positive roots for last reflections of~$w$, together with negative simple roots $-\alpha_s$ for each \mbox{$s\in S$} not appearing in $a$.
This set does not depend on the choice of reduced word for~$c$, because any two $c$-sorting words for~$w$ are related by commutations of simple generators.
One of the main results of~\cite{sortable} is the following theorem, which is an abbreviated form of \cite[Theorem~8.1]{sortable}.

\begin{theorem}
\label{cl}
The map $w\mapsto\cl_c(w)$ is a bijection from the set of $c$-sortable elements to the set of $c$-clusters.
\end{theorem}

\begin{example}
\label{B2cl}
In the case of $W=B_2$ and $c=s_0s_1$, $\cl_c(1)=\set{-\alpha_{s_0}, -\alpha_{s_1}}$ and the table below shows $\cl_c(w)$ for the other $c$-sortable elements~$w$.
\[\begin{array}{|c||c|c|c|c|c|}\hline
&&&&&\\[-7 pt]
w	&s_0		&s_0s_1		&s_0s_1s_0	&s_0s_1s_0s_1	&s_1\\[3 pt]\hline
&&&&&\\[-7 pt]
\cl_c(w)	& \alpha_{s_0}, -\alpha_{s_1}	& \alpha_{s_0}, \alpha_{s_0s_1s_0}	& \alpha_{s_1s_0s_1}, \alpha_{s_0s_1s_0} &\alpha_{s_1s_0s_1}, \alpha_{s_1}	&-\alpha_{s_0}, \alpha_{s_1}\\[3pt]\hline
\end{array}\]
\end{example}

The key result of this paper (Theorem~\ref{cl iso}) strengthens Theorem~\ref{cl} by asserting that $\cl_c$ induces a combinatorial isomorphism from~$\F_c$ to the $c$-cluster fan.
Via this combinatorial isomorphism, the $c$-Cambrian lattice induces a partial order on $c$-clusters which we call the {\em $c$-cluster lattice}.
The covering pairs of the $c$-cluster lattice are adjacent maximal cones of the $c$-cluster fan.
One moves down in the partial order by exchanging an almost positive root for another almost positive root which is ``closer'' to being a negative simple root, in a sense that is made precise in Section~\ref{lattice sec}.

\begin{figure}[ht]
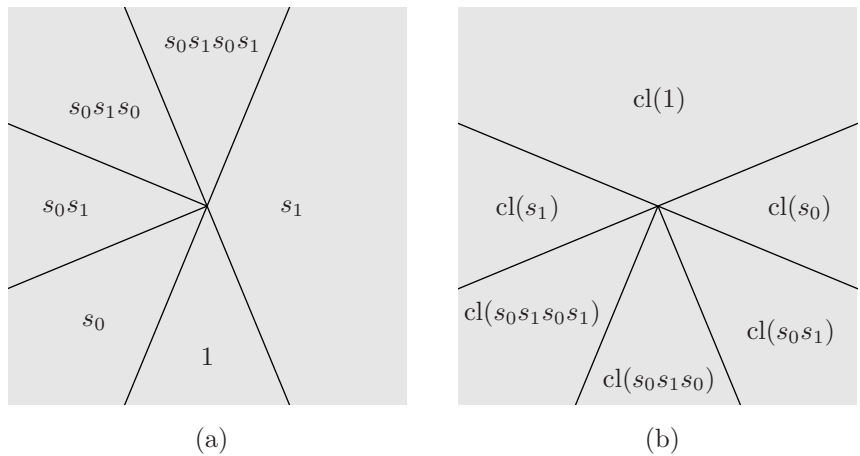

\centerline{
\begin{tabular}{cc}
\scalebox{1}{
\epsfbox{B2Fc.ps}
\begin{picture}(0,0)(75,-75)
\put(-6,-60){$1$}
\put(-51,-45){$s_0$}
\put(-66,-2){$s_0s_1$}
\put(-56,35){$s_0s_1s_0$}
\put(-20,60){$s_0s_1s_0s_1$}
\put(24,-2){$s_1$}
\end{picture}
}
&
\scalebox{1}{
\epsfbox{B2cclus.ps}
\begin{picture}(0,0)(75,-75)
\put(-13,38){$\cl(1)$}
\put(38,-3){$\cl(s_0)$}
\put(30,-50){$\cl(s_0s_1)$}
\put(-24.5,-68){$\cl(s_0s_1s_0)$}
\put(-65,-3){$\cl(s_1)$}
\put(-77,-43){$\cl(s_0s_1s_0s_1)$}
\end{picture}
}
\\[4 pt]
(a)&(b)
\end{tabular}
}
\caption{The $c$-Cambrian fan and the $c$-cluster fan}
\label{B2fan}
\end{figure}

\begin{example}
For $W=B_2$ with $c=s_0s_1$, the fan~$\F_c$ is shown in Figure~\ref{B2fan}.a, with maximal cones labeled by $c$-sortable elements.
Figure~\ref{B2fan}.b shows the $c$-cluster fan in the same coordinate system.
Each maximal cone in the $c$-cluster fan corresponds to the $c$-cluster composed of the extreme rays of the cone.
These maximal cones are labeled $\cl_c(w)$ (with the subscript~$c$ suppressed) for appropriate $c$-sortable elements~$w$.
The labeling of the rays is given in Figure~\ref{B2clusfan}. 
Observe that the obvious linear isomorphism between~$\F_c$ and the $c$-cluster fan is not induced by the bijection $\cl_c$.  
However, this linear isomorphism is an instance of a general result (Theorem~\ref{L iso}), which constructs, for special Coxeter elements called {\em bipartite} Coxeter elements, a linear isomorphism from the $c$-cluster fan to~$\F_c$.
\end{example}

The following simple lemma, which is \cite[Lemma~8.5]{sortable}, is a key ingredient in the proof (in~\cite{sortable}) of Theorem~\ref{cl} and in the results of this paper.
\begin{lemma}
\label{cl s}
Let~$s$ be initial in~$c$ and let~$w$ be $c$-sortable.
If $s\not\le w$ then $w\in W_{\br{s}}$ and $\cl_c(w)=\set{-\alpha_s}\cup\cl_{sc}(w)$.
If $s\le w$ then $\cl_c(w)=\sigma_s(\cl_{scs}(sw))$.
\end{lemma}

We have now surveyed the relevant background material. In Section~\ref{ray sec}, we undertake a detailed study of the rays of the Cambrian fan and prove the main technical lemmas underlying our main results.
In Section~\ref{proof sec}, we prove Theorem~\ref{cl iso}. 
We spend the remaining sections developing the further results described in Section~\ref{intro}.

\section{Rays in the $c$-Cambrian fan}
\label{ray sec}
In this section we prove the key lemmas which are used in the proof of Theorem~\ref{cl iso}.
These key lemmas, together with certain facts established in previous sections, can be loosely summarized as follows:
For~$s$ initial in~$c$, all of the objects relevant to Theorem~\ref{cl iso} are well-behaved under the operation of replacing~$c$ by $scs$ or, in some cases, replacing~$c$ by $sc$ and passing to the standard parabolic subgroup $W_{\br{s}}$. 
In particular, we define a map~$\zeta_s$ which describes how rays of the $c$-Cambrian fan transform under replacing~$c$ by $scs$.
We show that this map is compatible with the map~$\sigma_s$ on almost positive roots (defined in Section~\ref{cluster sec}). 

The results of this section rely (indirectly through results proved or quoted in Section~\ref{sort sec}) on nontrivial results from~\cite{sortable} and~\cite{sort_camb}.
We now specialize the description of the faces of $\F$ and~$\F_\Theta$ in the preceding section in order to describe the rays of~$\F_c$. 

Rays in the Coxeter fan $\F$ are in bijection with pairs $(w,J)$, where $w\in W$ and $J\subseteq S$ satisfy $|J|=n-1$ and $\ell(ws)>\ell(w)$ for every $s\in J$.
In particular~$w$ is either~$1$ or is a join-irreducible element of~$W\!$.
The correspondence is as follows:  
For any \mbox{$s\in S$}, let $\rho_s$ be the ray in the Coxeter fan which is fixed by $W_{\br{s}}$ and which is an extreme ray
of the region for~$1$. Note that~$\rho_s$ is usually not~$\alpha_s$. 
Given $(w,J)$, the corresponding ray is $w\rho_{s'}$, where $s'$ is the unique element of $S\setminus J$.
We write $\rho(w,J)$ for the ray associated to $(w,J)$.
Starting with a ray $\rho$ in the Coxeter fan, we recover $(w,J)$ as follows:
The elements of~$W$ whose regions contain $\rho$ form an interval in $W,$ and~$w$ is the minimal element of that interval.
The set $J$ is uniquely defined by specifying that the elements covering~$w$ in that interval are $\set{ws:s\in J}$. 

The following alternate description of $\rho(w,J)$ is also useful.
Namely, $\rho(w,J)$ is half of the line $I$ defined as the intersection of the hyperplanes associated to the reflections $\set{wsw^{-1}:s\in J}$.
Note that any reflecting hyperplane in~$\A$ either contains $I$ or intersects $I$ only at the origin.
For $w>1$, $\rho(w,J)$ is the half of $I$ consisting of points weakly separated from $D$ (the region for~$1$) by any hyperplane in~$\A$ which separates $wD$ from $D$.
For $w=1$, $\rho(w,J)$ is the half of $I$ which is contained in $D$. 

Proposition~\ref{faces of congruence} implies that rays in the $c$-Cambrian fan are the rays of the form $\rho(w,J)$ where~$w$ is $c$-antisortable.
By Corollary~\ref{num irr}, there are $|T|$ such pairs with $w\neq 1$.
There are also~$n$ such pairs with $w=1$, namely $(1,\br{r})$ for each $r\in S$.
We now proceed to define and then motivate a bijection~$\phi_c$ from rays of the $c$-Cambrian fan to almost positive roots. An example is given below (Example~\ref{vwJ}).

Given a ray $\rho(w,J)$ of the $c$-Cambrian fan, define $v=\pidown^c(w)$. 
In the case $w=1$ let \mbox{$\phi_c(\rho)=-\alpha_{s'}$} where $J=\br{s'}$. 
If $w\neq 1$ then~$w$ is join-irreducible, so by Proposition~\ref{irr bij},~$v$ is join-irreducible as well. 
Define $\phi_c(\rho)=\alpha_t$, where~$t$ is the unique cover reflection of~$v$.

We now motivate the definition of~$\phi_c$ by showing that it is forced on us by the requirement that $\phi_c^{-1}(\cl_c(v))$ be the set of rays of the cone $[v]_c$.  
If $w=1$ then $v=1$ and  $\cl_c(v)$ consists of all the negative simple roots for~$W\!$. 
In this case we must set $\phi_c(\rho)=-\alpha_s$ for some simple reflection \mbox{$s\in S$}. 
If instead we take $s\neq s'$ then $\phi_c^{-1}(-\alpha_s)=\rho_{s'}$ is a ray both of $[1]_c$ and $[s]_c$.
This is inconsistent with the fact that $-\alpha_s\not\in\cl_c(s)$.

To motivate the definition of~$\phi_c$ in the case $w\neq 1$, let~$v$ and~$t$ be as defined two paragraphs earlier.
Then any reduced word for~$v$ must end in the unique letter $r\in S$ such that $t=vrv^{-1}$.
In particular the $c$-sorting word for~$v$ must end in~$r$, so that~$t$ is the last reflection for~$r$ in~$v$.
The hyperplane~$H_t$ separates $[v]_c$ from the unique congruence class $[x]_c$ covered by $[v]_c$.
Since every region in $[x]_c$ is below the hyperplane~$H_t$, in particular~$t$ is not an inversion of $\pidown^c(x)$ and so $\alpha_t\not\in\cl_c(\pidown^c(x))$.
If $\cl_c$ is to induce a combinatorial isomorphism from the $c$-Cambrian fan to the $c$-cluster fan (which is simplicial) then $\cl_c(\pidown^c(x))$ and $\cl_c(v)$ should have $n-1$ roots in common and $\alpha_t$ should be the only root in $\cl_c(v)$ which is not in $\cl_c(\pidown^c(x))$.
Furthermore, the ray $\rho$ associated to $(w,J)$ should be the only ray of $[v]_c$ which is not a ray of $[x]_c$. Thus we are forced to map $\rho$ to $\alpha_t$.

To see that the map~$\phi_c$ is a bijection, note first that the~$n$ rays $\rho(1,J)$ map to the~$n$ negative simple roots.
The remaining rays are $\rho(w,J)$ where~$w$ is a $c$-antisortable join-irreducible, with a unique $J$ appearing for each such~$w$.
Thus Proposition~\ref{irr bij} shows that there is a unique~$v$ (equal to $\pidown^c(w)$) for each such pair and so 
by Corollary~\ref{irr cov},~$\phi_c$ is a bijection.

\begin{figure}
\centerline{\scalebox{.85}{
\epsfbox{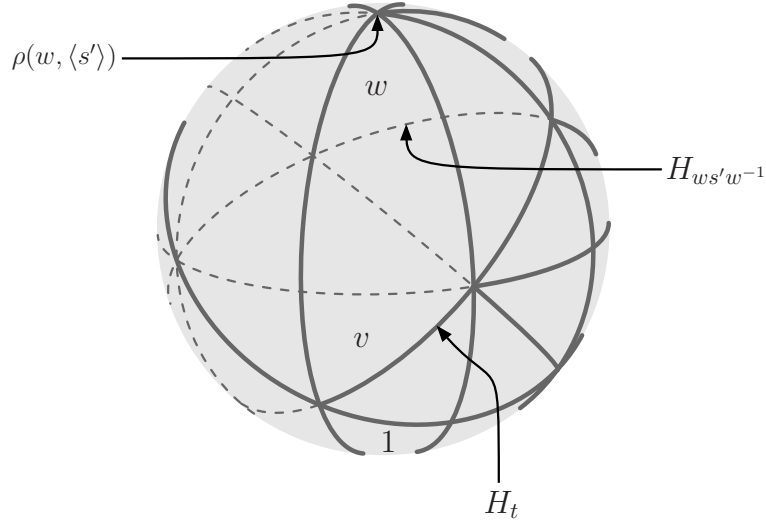}
\begin{picture}(0,0)(100,-135)
\put(-5,-104){\LARGE $1$}
\put(-17,-57){\LARGE $v$}
\put(-12,55){\LARGE $w$}
\put(-169,69){\Large $\rho(w,\br{s'})$}
\put(122,18){\LARGE $H_{ws'w^{-1}}$}
\put(41,-131){\LARGE $H_t$}
\end{picture}
}}
\caption{An illustration of the definition of~$\phi_c$}
\label{vwJfig}
\end{figure}

\begin{example}
\label{vwJ}
Figure~\ref{vwJfig} illustrates the definition of~$\phi_c$ in a typical instance of the case $w\neq 1$.
The solid lines show the intersection of the $c$-Cambrian fan with a unit hemisphere.
The dotted lines indicate how each maximal cone of the $c$-Cambrian fan is partitioned into maximal cones (regions) of the Coxeter fan.
Here~$s'$ is the unique element of $S\setminus J$ and $\phi_c(\rho(w,J))$ is the positive root associated to the hyperplane~$H_t$, the reflecting hyperplane for the unique cover reflection of $v$.
\end{example}

The following three lemmas constitute a recursive characterization of the rays of Cambrian fans.
They play a key role in the proof of Theorem~\ref{cl iso} in Section~\ref{proof sec}.

\begin{lemma}
\label{only ray}
If~$s$ is initial in~$c$ then $\rho_s$ is the only ray of the $c$-Cambrian fan which is below $H_s$ but not contained in $H_s$.
\end{lemma}
\begin{proof}
Suppose $\rho$ is a ray of the $c$-Cambrian fan which is on or below $H_s$ and suppose $\rho\neq\rho_s$.
We will show that $\rho$ is on $H_s$.
Let $\rho=\rho(w,J)$. 
The claim is easy when $(w,J)=(1,\br{r})$ for $r \neq s$ and we have excluded the case $(w,J)=(1,\br{s})$, so we may assume that $w \neq 1$.
Then since $\rho$ is on or below $H_s$,~$w$ is also below $H_s$, or in other words $s\not\le w$.
By Lemma~\ref{anti cov}, $sw\covers w$.
The ray $\rho$ is contained in the intersection of all hyperplanes separating~$w$ from an element covering~$w$, so in particular, $\rho$ is contained in $H_s$.
\end{proof}

\begin{lemma}
\label{para ray}
Let~$s$ be initial in~$c$. For any almost positive root $\alpha\in(\Phi_{\br{s}})_{\ge-1}$,
\[\phi_c^{-1}(\alpha)=\Proj^{-1}_{\br{s}}(\phi_{sc}^{-1}(\alpha)) \cap H_s,\]
where $H_s$ is the hyperplane associated to the reflection~$s$.
\end{lemma}

Recall that $\Proj_J$ is the linear projection $V(W) \onto V(W_J)$; in this case, the kernel of $\Proj_{\br{s}}$ is $\Span_{\RR} \rho_s$.
Here $\phi_c^{-1}(\alpha)$ is a ray of the $c$-Cambrian fan, and thus in particular a ray of the Coxeter fan for~$W\!$. The ray $\phi_{sc}^{-1}(\alpha)$ is the ray in the $sc$-Cambrian fan mapped to~$\alpha$ by $\phi_{sc}$.  Thus $\Proj^{-1}_{\br{s}}(\phi_{sc}^{-1}(\alpha))$ is two dimensional and $\Proj^{-1}_{\br{s}}(\phi_{sc}^{-1}(\alpha)) \cap H_{\br{s}}$ is a ray again.

Here is another description of Lemma~\ref{para ray}: The projection $\Proj_{\br{s}}$ restricts to an isomorphism from $H_s$ to $V_{\br{s}}$.
Identifying $V_{\br{s}}$ and $H_s$ by this isomorphism, Lemma~\ref{para ray} is the statement that $\phi_c^{-1}(\alpha)=\phi_{sc}^{-1}(\alpha)$.
However, this identification of $V_{\br{s}}$ with $H_s$ is not norm-preserving, and we will therefore not pursue this viewpoint.

%Another way to describe this result is to note that $\Proj_{\br{s}}$ gives an isomorphism from $H_s$ to $V_{\br{s}}$. Using this map to identify $V_{\br{s}}$ and $H_s$, Lemma~\ref{para ray} states that $\phi_c^{-1}(\alpha)$ is $\phi_{sc}^{-1}(\alpha)$. The reader should beware, however, that this identification of $V_{\br{s}}$ with $H_s$ is not norm-preserving, and we will therefore not pursue this viewpoint.

\begin{proof}[Proof of Lemma~\ref{para ray}]
Let $\rho(w,J)=\phi_c^{-1}(\alpha)$ and let $\rho(w',J')=\phi_{sc}^{-1}(\alpha)$.
Let $v=\pidown^c(w)$, let $v'=\pidown^{sc}(w')$  and let $\tilde{s}$ be the unique element of $S\setminus J$.
If $v=1$ then $\alpha=-\alpha_{\tilde{s}}$ and thus $v'=1$ as well, with $J'=J \setminus \set{s}$. 
Now $\phi_c^{-1}(\alpha)=\rho_{\tilde{s}}$, where $\rho_{\tilde{s}}$ is a ray in $V(W)$, and $\phi_{sc}^{-1}(\alpha)$ is also $\rho_{\tilde{s}}$, interpreted as a ray in $V(W_{\br{s}})$. 
The desired conclusion now follows from the definition of $\rho_{\tilde{s}}$ in each case.

If $v>1$ then~$\alpha$ is not a negative simple root, so $v'>1$ as well.
Thus $v'$ is join-irreducible, and since $v'$ is $sc$-sortable, it is in particular $c$-sortable.
For both~$v$ and $v'$, the positive root associated to the unique cover reflection is~$\alpha$.
We have (by Corollary~\ref{irr cov}) $v=v'$ and thus $w=\piup_c(v)$ and $w'=\piup_{sc}(v)$.  
The ray $\rho(w,J)$ is half of the line $I$ defined as the intersection of the hyperplanes $\set{H_{wrw^{-1}}(V(W)):r\in J}$ and the ray $\rho(w',J')$ is half of the line $I'$ defined as the intersection of the hyperplanes $\set{H_{w'r(w')^{-1}}(V(W_{\br{s}})):r\in J'}$. Here we have written $H_t(V(W))$ for the hyperplane in $V(W)$ associated to the reflection~$t$. We have $\Proj_{\br{s}}^{-1}(H_t(V(W_{\br{s}}))) = H_t(V(W))$.

We next show that $I=\Proj^{-1}_{\br{s}}(I') \cap H_s$ by proving the stronger statement 
\[\set{wrw^{-1}:r\in J}=\set{w'r(w')^{-1}:r\in J'}\cup\set{s}.\]
Recall that the set $\set{wrw^{-1}:r\in J}$ is the set of all reflections~$t$ such that~$H_t$ separates~$w$ from an adjacent region covering~$w$ in weak order. By the dual of Lemma~\ref{separators}, these are the reflections~$t$ such that~$H_t$ separates $[w]_c=[w']_c$ from an adjacent region covering it in~$\F_c$. Similarly, $\set{w'r(w')^{-1}:r\in J'}$ is the set of reflections $t \in W_{\br{s}}$ such that~$H_t$ separates $[w']_{sc}$ from an adjacent region covering it in $\F_{sc}$. But the $c$-Cambrian congruence restricted to $W_{\br{s}}$ is simply the $sc$-Cambrian congruence so any region covering $[w']_{sc}$ in $W/{\Theta_{sc}}$ also corresponds to a cover of $[w']_c$ in $W/{\Theta_c}$ with a separating hyperplane corresponding to the same reflection. Thus,  $\set{w'r(w')^{-1}:r\in J'}\subseteq\set{wrw^{-1}:r\in J}$.

Now, the set $\set{wrw^{-1}:r\in J}$ contains exactly one additional element not contained in $\set{w'r(w')^{-1}:r\in J'}$. Since $s\not\le v$, Lemma~\ref{no crossing} says that $s\not\le w$.
Thus Lemma~\ref{anti cov} says that $sw$ covers~$w$, so that $s\in\set{wrw^{-1}:r\in J}$, completing the proof that $\set{wrw^{-1}:r\in J}=\set{w'r(w')^{-1}:r\in J'}\cup\set{s}$.

Having established that $\phi_c^{-1}(\alpha)$ and $\Proj^{-1}_{\br{s}}(\phi_{sc}^{-1}(\alpha)) \cap H_s$ each constitute half of the line $I$, it remains to show that they are same half of $I$.
Since $v'>1$, there is a hyperplane~$H$ in $\A_{\br{s}}$ which separates $w'D$ from $D$, and since $w\ge w'$,~$H$ separates $wD$ from $D$ as well.
Thus both $\phi_c^{-1}(\alpha)$ and $\Proj^{-1}_{\br{s}}(\phi_{sc}^{-1}(\alpha)) \cap H_s$ are the half of $I$ weakly separated from $D$ by~$H$.
\end{proof}

For a ray $\rho$ of the Cambrian fan, define 
\[\zeta_s(\rho)=\left\lbrace\begin{array}{ll}
s\rho&\mbox{if }\rho\neq\rho_s,\mbox{ or}\\
-\rho&\mbox{if }\rho=\rho_s.\end{array}\right.\]
Since the Coxeter fan is preserved by the action of~$s$ and by the antipodal map, $\zeta_s(\rho)$ is a ray of the Coxeter fan.\footnote{One should note that the inverse map $\zeta_s^{-1}$ from rays of the $scs$-Cambrian fan to rays of the $c$-Cambrian fan is \emph{not} given by the same formula.}
The following lemma states that $\zeta_s(\rho)$ is a ray of the $scs$-Cambrian fan and establishes the compatibility of~$\zeta_s$ with~$\phi_c$, $\sigma_s$ and~$\phi_{scs}$.
The proof is not difficult, but it has many cases.  

\begin{lemma}
\label{zeta}
Let~$s$ be an initial letter of~$c$ and let $\rho$ be a ray in the $c$-Cambrian fan.
Then $\zeta_s(\rho)$ is a ray in the $scs$-Cambrian fan and $\phi_{scs}(\zeta_s(\rho))=\sigma_s(\phi_c(\rho))$. 
\end{lemma}

\begin{proof}
Let~$s$,~$c$ and $\rho$ be as in the statement of the lemma and let $\rho=\rho(w,J)$. Further, let $s'$ be such that $J=\br{s'}$ and let $v=\pidown^c(w)$.

\textbf{Case 1:} $w=1$. This case splits into two subcases:

\textbf{Case 1a:}  $(w,J)=(1,\br{s})$. This is the exceptional case in the definition of~$\zeta_s$, where  $\rho=\rho_s$. In this case $\phi_c(\rho)=-\alpha_s$, so $\sigma_s(\phi_c(\rho))=\alpha_s$.
Since $\zeta_s(\rho)=-\rho_s$, we need to show that $-\rho_s$ is a ray in the $scs$-Cambrian fan and that $\phi_{scs}(-\rho_s)=\alpha_s$.
Let $\rho'=\rho(w',J')$ be the unique ray of the $scs$-Cambrian fan with $\phi_{scs}(\rho')=\alpha_s$.
Since~$s$ is the unique $scs$-sortable join-irreducible whose associated reflection is~$s$, we must have $\pidown^{scs}(w')=s$  so that $\piup_{scs}(s)=w'$.
Thus $w'=w_0\cdot\left((w_0)_{\br{w_0sw_0}}\right)$ by Lemma~\ref{piup s}. The ascents of $w'$ are $J'=\br{w_0sw_0}$, so $\rho'=w'\rho_{w_0sw_0}=w_0\rho_{w_0sw_0}$, the latter equality holding because $(w_0)_{\br{w_0sw_0}}$ fixes $\rho_{w_0sw_0}$. But $w_0\rho_{w_0sw_0}=-\rho_s$.

\textbf{Case 1b:}  $(w,J)=(1,\br{r})$ for $r\neq s$.
In this case $\rho=\rho_r$ and since $\rho_r$ is on the reflecting hyperplane for~$s$, $\zeta_s(\rho)=s\rho=\rho$.
Since~$1$ is also $scs$-sortable and $J$ is a set of $n-1$ ascents of~$1$, $\rho$ is also a ray of the $scs$-Cambrian fan with $\phi_{scs}(\rho)=-\alpha_r$.
Also, $\sigma_s(\phi_c(\rho))=\sigma_s(-\alpha_r)=-\alpha_r$.

\textbf{Case 2:}  $s\le w$. By Lemmas~\ref{s ji} and~\ref{scs dual}, $sw$ is $scs$-antisortable and $sw$ is either~$1$ or join-irreducible.
Furthermore, $J$ is a set of elements which lengthen not only~$w$ but also $sw$ on the right.
We consider two subcases, depending on whether or not $w=s$.
Notice that since~$s$ is initial in~$c$, the sole element in the $c$-Cambrian equivalence class of $s$ is $s$ itself.
Thus $w=s$ if and only if $v=s$.

\textbf{Case 2a:}  $w=s$.  In this case $(sw,J)=(1,\br{s})$, with associated ray $\rho_s=\zeta_s(\rho)$.
Also $v=s$ and $\rho=s\rho_s$.
So we have $\sigma_s(\phi_c(\rho))=\sigma_s(\alpha_s)=-\alpha_s=\phi_{scs}(\rho_s)=\phi_{scs}(\zeta_s(\rho))$.

\textbf{Case 2b:}  $w\neq s$. In this case $\phi_c(\rho)$ is a positive root $\alpha_t$ for some reflection $t\neq s$ and thus $\sigma_s(\phi_c(\rho))$ is $\alpha_{sts}$.
The ray $\rho(sw,J)$ of the $scs$-Cambrian fan is $\rho':=sw\rho_{s'}=\zeta_s(\rho)$.
Also, $\pidown^c(w)=s\cdot\pidown^{scs}(sw)$, so that $\pidown^{scs}(sw)=s\cdot\pidown^c(w)=sv$.
The reflection~$t$ is the unique cover reflection of~$v$ and, by Lemma~\ref{s ji}, the element $sv$ is join-irreducible with unique cover reflection $sts$.
In particular $\phi_{scs}(\rho')=\alpha_{sts}$. This concludes the proof for the case $s\le w$.

\textbf{Case 3:}  $s\not\le w$ and $w \neq 1$. 
By Lemma~\ref{anti cov},~$w$ is $scs$-antisortable and $sw=s \join w$.
Since~$w$ is $scs$-antisortable, the pair $(w,J)$ defines $\rho$ not only as a ray of the $c$-Cambrian fan but also as a ray of the $scs$-Cambrian fan.
Since $w\neq 1$ but~$w$ is below $H_s$, $\rho$ is contained in $H_s$  by Lemma~\ref{only ray} so $\rho=\zeta_s(\rho)$.

Let $v'=\pidown^{scs}(w)$.
By definition, $\pidown^c(sw)=s\cdot\pidown^{scs}(w)=sv'$.
Since $\pidown^c$ is a lattice homomorphism, $\pidown^c(sw)=\pidown^c(s\join w)=s\join\pidown^c(w)=s\join v$.
Thus $v'=s\cdot(s\join v)$.
Because $s\not\le w$,~$v$ is in $W_{\br{s}}$, so by Lemma~\ref{cov w br s}, the set of cover reflections of $s\join v$ is $\set{s,t}$, where~$t$ is the unique cover reflection of~$v$.
Let~$t'$ be the unique cover reflection of $v'$; by definition, $\phi_{scs}(w')=\alpha_{t'}$.
Since the interval $[1,sw_0]$ is isomorphic to the interval $[s,w_0]$ by the map $x\mapsto sx$, the reflection $st's$ is a cover reflection of $sv'=s\join v$.
But $v'\not\ge s$, so $t'\neq s$ and thus $st's=t$.
Therefore $\phi_{scs}(\rho)=\alpha_{sts}=\sigma_s(\alpha_t)=\sigma_s(\phi_c(\rho))$.
\end{proof}

\section{Proof of the combinatorial isomorphism}
\label{proof sec}
In this section we prove the main theorem, Theorem~\ref{cl iso}, which states that, for~$W$ finite, the map $\cl_c$ induces a combinatorial isomorphism from the $c$-Cambrian fan to the $c$-cluster fan. 
We also discuss some first consequences of Theorem~\ref{cl iso}.

The $c$-cluster fan is simplicial.  
That is, each of its maximal faces is the positive linear span of a collection of linearly independent vectors.
Specifically, this collection of vectors is a $c$-cluster of almost positive roots.
In contrast, we do not even know that the $c$-Cambrian fan is simplicial.
However, we know that the maximal cones of the $c$-Cambrian fan are, by definition, unions of regions of the Coxeter arrangement.  (See Section~\ref{camb fan sec}.)
Specifically, each maximal cone of the $c$-Cambrian fan is the union over a fiber of the map $\pidown^c$.
Showing that the $c$-Cambrian fan is simplicial means showing the following:
For each fiber of $\pidown^c$, there is a collection $E$ of~$n$ rays of the $c$-Cambrian fan such that  a given region is a member of the fiber if and only if that region is  contained in the positive linear span of $E$.

The stronger statement that $\cl_c$ induces a combinatorial isomorphism between the $c$-Cambrian fan and the $c$-cluster fan is equivalent to the additional condition that there is a bijection $\phi$ between the rays of the $c$-Cambrian fan and the rays of the $c$-cluster fan such that the collection $E$ of rays used to determine membership in $(\pidown^c)^{-1}(x)$ obeys $\phi(E)=\cl_c(x)$.

As was shown in Section~\ref{ray sec}, the map~$\phi_c$ is a bijection from rays of the $c$-Cambrian fan to almost positive roots---that is, to rays of the $c$-cluster fan.
Thus the proof of Theorem~\ref{cl iso} is completed by Proposition~\ref{span} below. 
Recall that, if $W_J$ is a standard parabolic subgroup, then the dominant chamber of $\A(W_J)$ is denoted by $D_J$.

%We remind the reader that there is a region $D$ of the hyperplane arrangement for $W$, called the dominant chamber and defined by $\bigcap_{s \in S} \{ v : \langle v,\alpha_s \rangle \geq 0 \}$ and that, to each $w \in W$, we associate the region $wD$ of the hyperplane arrangement. If $W_J$ is a standard parabolic subgroup, we denote the dominant chamber of $W_J$ by $D_J$.

\begin{prop} \label{span}
Let~$x$ be $c$-sortable.
Then the following are equivalent for any $w\in W$.
\begin{enumerate}
\item[(i) ] $\pidown^c(w)=x$.
\item[(ii) ] The interior of the region $wD$ intersects the positive span of $\phi_c^{-1}(\cl_c(x))$.
\item[(iii) ] The region $wD$ is contained in the positive span of $\phi_c^{-1}(\cl_c(x))$.
\end{enumerate}
\end{prop}
Here, since $\cl_c(x)$ is a cluster of almost positive roots, $\phi_c^{-1}(\cl_c(x))$ represents the set of rays obtained by applying $\phi_c^{-1}$ to each member of the cluster.

\begin{proof}
%For convenience in this proof, we tacitly identify elements~$w$ of~$W$ with regions $wD$ of the Coxeter arrangement for~$W\!$.
The fact that (iii) implies (ii) is trivial.
We prove that (i) implies (iii) and that (ii) implies (i) by induction on the length of~$w$ and the rank of~$W\!$.
Let~$s$ be initial in~$c$.  
For each implication we will consider two cases: $s\not\le w$ and $s\le w$.

First, assume (i).
If $s\not\le w$ then $\pidown^{sc}(w_{\br{s}})=x$, so that in particular $x\in W_{\br{s}}$.
By Lemma~\ref{cl s}, $\cl_c(x)=\cl_{sc}(x)\cup\set{-\alpha_s}$.
By Lemma~\ref{para ray}, each ray in $\phi_c^{-1}(\cl_{sc}(x))$ is obtained from the corresponding ray $\rho$ in $\phi_{sc}^{-1}(\cl_{sc}(x))$ by intersecting $\Proj^{-1}_{\br{s}}(\rho)$ with the hyperplane $H_s$.  
Since $\phi_c^{-1}(-\alpha_s)=\rho_s$ is the half of the intersection of the hyperplanes in $\A_{\br{s}}$ which is below the hyperplane $H_s$, the positive span of $\phi_c^{-1}(\cl_c(x))$ is the part of the positive span of $\Proj^{-1}_{\br{s}}(\phi_{sc}^{-1}(\cl_{sc}(x)))$ which is below the hyperplane $H_s$. 
Now,~$wD$ is contained in $\Proj^{-1}_{\br{s}}(w_{\br{s}} D_{\br{s}})$ (Proposition~\ref{projection}) which is, by induction on rank, contained in the positive span of $\Proj^{-1}_{\br{s}}(\phi_{sc}^{-1}(\cl_{sc}(x)))$. 
Since $s\not\le w$,~$wD$ is below the hyperplane $H_s$. Thus we see that~$wD$ is in the positive span of $\phi_c^{-1}(\cl_c(x))$.  

If $s\le w$ then (i) implies $\pidown^{scs}(sw)=sx$.
By Lemma~\ref{no crossing}, $s\le x$, and thus $x\not\in W_{\br{s}}$, so that in particular $\rho_s\not\in\phi_c^{-1}(\cl_c(x))$.
By induction on length, $swD$ is completely in the positive span of $\phi_{scs}^{-1}(\cl_{scs}(sx))$.
By Lemma~\ref{cl s}, $\phi_{scs}^{-1}(\cl_{scs}(sx))$ equals $\phi_{scs}^{-1}(\sigma_s(\cl_c(x)))$, which by Lemma~\ref{zeta} equals $\zeta_s(\phi_c^{-1}(\cl_c(x)))=s\phi_c^{-1}(\cl_c(x))$.
Since $swD$ is completely in the positive span of $s\phi_c^{-1}(\cl_c(x))$,~$wD$ is completely in the positive span of $\phi_c^{-1}(\cl_c(x))$.

Now suppose (ii).
If $s\not\le w$ then every point in the interior of~$wD$ is strictly below $H_s$.
By Lemma~\ref{only ray}, $\phi_c^{-1}(\cl_c(x))$ must contain the ray $\rho_s$, so that $\cl_c(x)$ contains $-\alpha_s$.
In particular, $x\in W_{\br{s}}$, and furthermore by Lemma~\ref{cl s}, $\cl_c(x)=\cl_{sc}(x)\cup\set{-\alpha_s}$.
By Proposition~\ref{projection},  $\Proj_{\br{s}}(wD) \subseteq w_{\br{s}} D_{\br{s}}$ and the interior of $wD$ is taken into the interior of $w_{\br{s}} D_{\br{s}}$ by $\Proj_{\br{s}}$ (by considerations of dimension).  
Since $\rho_s$ is in the kernel of $\Proj_{\br{s}}$, the interior of $w_{\br{s}} D_{\br{s}}$ intersects the positive span of $\phi_{sc}^{-1}(\cl_{sc}(x))$. By induction on rank, $\pidown^{sc}(w_{\br{s}})=x$ and thus $\pidown^c(w)=x$.

If $s\le w$ then we claim that $s\le x$.
Supposing to the contrary that $s\not\le x$, by Lemma~\ref{cl s}, $\cl_c(x)=\set{-\alpha_s}\cup\cl_{sc}(x)$.
Thus by Lemma~\ref{para ray}, $\phi_c^{-1}(\cl_c(x))$ consists of rays which are weakly below $H_s$.
But the interior of~$wD$ is strictly above $H_s$, contradicting the supposition that (ii) holds.
This contradiction proves the claim that $s\le x$.
In particular, $\rho_s$ is not in $\phi_c^{-1}(\cl_c(x))$, so that $\zeta_s(\phi_c^{-1}(\cl_c(x)))=s\phi_c^{-1}(\cl_c(x))$.
Thus the interior of $swD$ meets the positive span of $\zeta_s(\phi_c^{-1}(\cl_c(x)))$, which by Lemma~\ref{zeta} equals $\phi_{scs}^{-1}(\sigma_s(\cl_c(x)))$.
Since $s\le x$, Lemma~\ref{cl s} says that the latter is $\phi_{scs}^{-1}(\cl_{scs}(sx))$.
By induction on length, $\pidown^{scs}(sw)=sx$, so that \mbox{$\pidown^c(w)=x$}.
\end{proof}

This completes the proof of Theorem~\ref{cl iso}. 
In fact, we have proven the following more detailed version of Theorem~\ref{cl iso}.
\begin{theorem} \label{detailed cl iso}
The $c$-Cambrian fan~$\F_c$ is simplicial and the bijection~$\phi_c$ between the rays of the Cambrian fan and the almost positive roots induces a combinatorial isomorphism of fans between the $c$-Cambrian fan and the $c$-cluster fan. Under this isomorphism, the maximal cone $[w]_c$ is taken to the cluster $\cl_c(w)$.
\end{theorem}

If~$w$ is $c$-sortable and $x\covered w$ then the maximal cones $[w]_c$ and $[\pidown^c(x)]_c$ intersect in a facet of dimension $n-1$ and so have exactly $n-1$ rays in common.
Thus Theorem~\ref{cl iso} has the following corollary.
\begin{cor}
\label{cl compat}
Let~$w$ be $c$-sortable and let $x\covered w$.
Then the $c$-clusters $\cl_c(w)$ and $\cl_c(\pidown^c(x))$ have exactly $n-1$ almost positive roots in common.
\end{cor}

In Section~\ref{cluster sec}, we noted that the action of~$\sigma_s$ on almost positive roots induces a combinatorial isomorphism between the $c$-cluster fan and the $scs$-cluster fan.
Thus, by Theorem~\ref{cl iso}, the map $\phi^{-1}_{scs}\circ\sigma_s\circ\phi_c$ induces a combinatorial isomorphism between~$\F_c$ and~$\F_{scs}$.
But Lemma~\ref{zeta} implies that $\phi^{-1}_{scs}\circ\sigma_s\circ\phi_c$ coincides with $\zeta_s$.
Since the combinatorial isomorphism is determined by its action on rays, we have the following.
\begin{prop} \label{zeta and sigma}
The action of~$\zeta_s$ on the rays of the $c$-Cambrian fan~$\F_c$ induces a combinatorial isomorphism between~$\F_c$ and~$\F_{scs}$.
\end{prop}
In particular~$\F_c$ and~$\F_{scs}$ are related by a piecewise-linear map that is only a slight deformation of the linear map~$s$.
On and above $H_s$ the map agrees with~$s$.
Below $H_s$ the maps agrees with the linear map that fixes $H_s$ and takes $\rho_s$ to $-\rho_s$.

We now describe the isomorphism between~$\F_c$ and~$\F_{scs}$ directly in terms of sortable elements (Cf.\ \cite[Remark~3.8]{sort_camb}). For~$s$ initial in~$c$, define a map~$Z_s$ from the set of $c$-sortable elements of~$W$ to the set of $scs$-sortable elements by:
\[Z_s(w) =\left\lbrace\begin{array}{ll}
sw & \mbox{if } s\le w,\mbox{ or}  \\
s \join w &\mbox{if } s\not\le w.
\end{array}\right.\]
We now check that~$Z_s$ maps $c$-sortable elements to $scs$-sortable elements: If $s\le w$, this is Lemma~\ref{scs}. If $s\not\le w$ then Lemma~\ref{sc} states that $w \in W_{\br{s}}$ and~$w$ is $sc$-sortable. The $sc$-sorting word for~$w$ is identically equal to the $scs$-sorting word for~$w$, so that~$w$ is $scs$-sortable. The set of $scs$-sortable elements forms a sublattice of~$W$ (Theorem~\ref{sublattice}), so $s \join w$ is also $scs$-sortable. 

The inverse of~$Z_s$ is
\[Z^{-1}_s(w) =\left\lbrace\begin{array}{ll}
w_{\br{s}} & \mbox{if } s\le w,\mbox{ or}  \\
sw &\mbox{if } s\not\le w.
\end{array}\right.\]
There are only two nontrivial assertions in the statement that this map is indeed the inverse of~$Z_s$:
first that any $c$-sortable element~$w$ with $s\not\le w$ obeys the condition $(s\join w)_{\br{s}}=w$; and second that an $scs$-sortable element~$w$ with $s\le w$ obeys the condition $(w_{\br{s}}\join s)=w$.
Recall that $x\mapsto x_{\br{s}}$ is a lattice homomorphism, so that $(s\join w)_{\br{s}}=s_{\br{s}}\join w_{\br{s}}=w_{\br{s}}$.
Thus the first assertion follows from Lemma~\ref{sc}.
The second assertion is exactly \cite[Lemma~2.10]{sort_camb}.

The following lemma states that $[w]_c \mapsto [Z_s(w)]_{scs}$ is the isomorphism between~$\F_c$ and~$\F_{scs}$ corresponding to the isomorphism~$\sigma_s$ of cluster fans.

\begin{lemma}
\label{BigZ}
For a $c$-sortable element~$w$, $\cl_{scs}(Z_s(w))=\sigma_s\cl_c(w)$.
\end{lemma}
\begin{proof}
If $s\le w$ then Lemma~\ref{cl s} is the desired statement.

If $s\not\le w$ then the desired equality is $\cl_{scs}(s\join w)=\sigma_s\cl_c(w)$.
Since $\phi_{scs}$ is a bijection, this is equivalent to checking that $\phi_{scs}^{-1}\cl_{scs}(s\join w)=\phi_{scs}^{-1}\sigma_s\cl_c(w)$, which can be rewritten, using Lemma~\ref{zeta}, as $\phi_{scs}^{-1}\cl_{scs}(s\join w)=\zeta_s\phi_c^{-1}\cl_c(w)$.
In other words, the requirement is that the rays $\rho_1,\ldots,\rho_n$ of the $c$-Cambrian cone $[w]_c$ are mapped by~$\zeta_s$ to the rays of the $scs$-Cambrian cone $[s\join w]_{scs}$. By Lemma~\ref{no crossing} all the $\rho_i$ are below $H_s$. By Lemma~\ref{only ray}, all of the $\rho_i$ are in $H_s$ except for possibly one, which is $\rho_s$. We know that $\rho_1$, \dots, $\rho_n$ are linearly independent, so one of the $\rho_i$ must be $\rho_s$; without loss of generality let $\rho_n=\rho_s$. Then $\zeta_s(\rho_n)=-\rho_s$ and $\zeta_s(\rho_i)=\rho_i$ for $i<n$. Now, for any $u \in W,$ $u$ is in the positive span of $-\rho_s$ and \mbox{$\rho_1,\ldots,\rho_{n-1}$} if and only if the following conditions hold: $u \geq s$ and $u_{\br{s}}$, considered as a region of $V(W_{\br{s}})$, is in the positive span of \mbox{$\rho_1,\ldots,\rho_{n-1}$}. 

We have $s \join w \geq s$ and also $(s \join w)_{\br{s}}=s_{\br{s}} \join w_{\br{s}}=w_{\br{s}}=w$. 
Our hypothesis is that~$w$, when considered as a region of $V(W)$, is in the positive span of \mbox{$\rho_1,\ldots,\rho_n$}.
This implies that~$w$  considered as a region of $V(W_{\br{s}})$ is in the positive span of \mbox{$\rho_1,\ldots,\rho_{n-1}$}. 
So we conclude that $s \join w$ is in the positive span of the $\zeta_s(\rho_i)$ as desired. 
\end{proof}

\begin{figure}
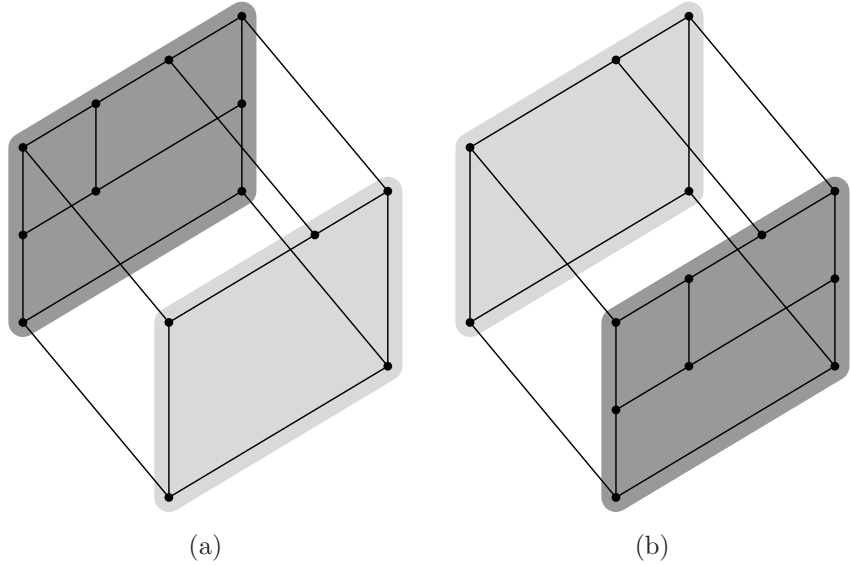

\centerline{
\begin{tabular}{ccc}
\scalebox{.55}{\epsfbox{cambA3abc.ps}}
&&
\scalebox{.55}{\epsfbox{cambA3bca.ps}}
\\[4 pt]
(a)&&(b)
\end{tabular}
}
\caption{The $c$- and $scs$-Cambrian lattices}
\label{A3Zfig}
\end{figure}

\begin{example}
\label{A3Z}
The map~$Z_s$ is perhaps more easily visualized as a map from the $c$-Cambrian lattice to the $scs$-Cambrian lattice.
Figure~\ref{A3Zfig}.a shows the $s_1s_2s_3$-Cambrian lattice for~$W$ of type $A_3$; this is also the lattice depicted in Figure~\ref{A3cambfig}.b. 
The light gray shading indicates (congruence classes of) $s_1s_2s_3$-sortable elements not above~$s_1$, while dark gray shading indicates $s_1s_2s_3$-sortable elements above~$s_1$.
Figure~\ref{A3Zfig}.b shows the $s_2s_3s_1$-Cambrian lattice for the same~$W\!$.
The map $Z_{s_1}$ takes the $s_1s_2s_3$-sortable elements not above~$s_1$ to the $s_2s_3s_1$-sortable elements above~$s_1$, which are shaded light gray in Figure~\ref{A3Zfig}.b.
The $s_1s_2s_3$-sortable elements above~$s_1$ are taken to $s_2s_3s_1$-sortable elements not above~$s_1$, shaded dark gray in Figure~\ref{A3Zfig}.b.
Notice that $Z_{s_1}$ restricted to light-shaded elements in Figure~\ref{A3Zfig}.a is a poset isomorphism to light-shaded elements in Figure~\ref{A3Zfig}.b, and similarly for dark-shaded elements.
\end{example}

\section{The cluster lattice}
\label{lattice sec}
We have now established in great detail the combinatorial isomorphism between the $c$-Cambrian fan and the $c$-cluster fan. 
The maximal cones of the $c$-Cambrian fan are partially ordered by the Cambrian lattice $W/\Theta_c$, so we obtain an induced poset (in fact lattice) structure on the set of $c$-clusters. 
In this section we will apply our results to describe this poset directly in terms of cluster combinatorics.

The {\em exchange graph} on $c$-clusters is the adjacency graph on maximal cones of the $c$-cluster complex.
In other words, the vertices are the $c$-clusters, with an edge between~$C$ and~$C'$ if and only if $| C \cap C' |=n-1$.
The exchange graph is isomorphic to the $1$-skeleton of the (simple) generalized associahedron for~$W$ as defined in~\cite{ga}.
We have the following corollary of Theorem~\ref{cl iso}.

\begin{cor} \label{Hasse}
The undirected Hasse diagram of the $c$-Cambrian lattice $W/\Theta_c$ is isomorphic to the exchange graph on $c$-clusters and hence isomorphic to the $1$-skeleton of the generalized associahedron for~$W\!$.
\end{cor}
\begin{proof}
Proposition~\ref{ItsACover} implies that the Hasse diagram of $W/\Theta_c$ is the adjacency graph of maximal cones of~$\F_c$, which, according to Theorem~\ref{cl iso}, is mapped by $\cl_c$ to the exchange graph.
\end{proof}

In light of Corollary~\ref{Hasse}, to describe the poset induced on $c$-clusters by the $c$-Cambrian lattice, it is sufficient to give the correct orientation of the exchange graph.
Two almost positive roots~$\alpha$ and $\alpha'$ are said to be {\em $c$-exchangeable} if they are distinct and if there is a set $B\subseteq(\Pge\setminus\set{\alpha,\alpha'})$ such that both $B \cup \set{\alpha}$ and $B \cup \set{\alpha'}$ are $c$-clusters.
Note that $c$-exchangeable roots~$\alpha$ and $\alpha'$ are never $c$-compatible.
We will show that the correct orientation of an edge $B \cup \set{\alpha} \hbox{\ \textemdash\ } B \cup \set{\alpha'}$ depends only on~$c$ and the pair $(\alpha, \alpha')$ of $c$-exchangeable roots.
Specifically, the orientation is given by comparing $\alpha $ and $\alpha'$ using a function $R_c$ from almost positive roots to natural numbers which we now proceed to define.

Recall from the introduction the involution $\sigma_s : \Pge \to \Pge$ for each \mbox{$s\in S$}. 
Let $s_1 \cdots s_n$ be a reduced word for~$c$ and define $\sigma_c=\sigma_{s_1} \cdots \sigma_{s_n}$. 
Any two reduced words for~$c$ differ only by interchanging commuting reflections and $\sigma_{s} \sigma_{t}=\sigma_{t} \sigma_{s}$ whenever $st=ts$, so $\sigma_c$ is a well defined permutation of $\Pge$. Note that \mbox{$(\sigma_c)^{-1}=\sigma_{c^{-1}}$}.

\begin{prop} \label{BigR}
For any almost positive root~$\alpha$ and any Coxeter element~$c$, there exists a nonnegative integer $R$ such that $\sigma_c^{-R}(\alpha)$ is a negative simple root. 
\end{prop}

Proposition~\ref{BigR} will be proved later in the section. We write $R_c(\alpha)$ for the smallest such $R$. Assuming the proposition, we define the {\em $c$-cluster lattice} $\Clust_c$ to be the partial order on $c$-clusters whose cover relations are $B\cup\set{\alpha}\covered B\cup\set{\alpha'}$ if and only if $R_c(\alpha)<R_c(\alpha')$.
(Cf.\ \cite[Section~8]{cambrian}).
It is not obvious from this definition that these relations are in fact cover relations of the partial order they generate.
However, in light of the preceding discussion, the following proposition implies that the relations above are in fact cover relations.

\begin{prop}
\label{cl covers}
Suppose $\alpha$ and $\alpha'$ are $c$-exchangeable almost positive roots and let $B$ be a subset of $(\Pge\setminus\set{\alpha,\alpha'})$ such that $B\cup\set{\alpha}$ and $B\cup\set{\alpha'}$ are $c$-clusters.
Then $[cl^{-1}_c(B\cup\set{\alpha})]_c$ is covered by $[cl^{-1}_c(B\cup\set{\alpha'})]_c$ in the $c$-Cambrian lattice if and only if $R_c(\alpha)<R_c(\alpha')$.
\end{prop}
Notice that the case $R_c(\alpha)=R_c(\alpha')$ is impossible for~$\alpha$ and $\alpha'$ as in the proposition.
This is because if $R_c(\alpha)=R_c(\alpha')$, one can iterate the definition of $c$-compatibility (see Section~\ref{cluster sec}) to show that $\alpha \cm_c \alpha'$,  contradicting the fact that~$\alpha$ and $\alpha'$ are $c$-exchangeable.
Proposition~\ref{cl covers} and Corollary~\ref{Hasse} immediately imply the main theorem of this section, which is a generalization of the second statement of \cite[Conjecture~1.4]{cambrian}.

\begin{theorem}
\label{cl poset iso}
The map $\cl_c$ is an isomorphism of lattices from the $c$-Cambrian lattice $W/\Theta_c$ to the $c$-cluster lattice.
\end{theorem}

Before proving Propositions~\ref{BigR} and~\ref{cl covers}, we mention some consequences of Theorem~\ref{cl poset iso}.
The following corollary is immediate from Theorem~\ref{cl iso}, Theorem~\ref{cl poset iso} and the properties of~$\F_c$ and $W/\Theta_c$ listed in Section~\ref{camb fan sec}.

\begin{cor}
\label{consequences}
The $c$-cluster lattice $\Clust_c$ has the following properties.
\begin{enumerate}
\item[(i) ] $\Clust_c$ is a lattice.
\item[(ii) ]  Any linear extension of $\Clust_c$ is a shelling order of the $c$-cluster complex.
\item[(iii) ] For any set $F$ of $c$-compatible almost positive roots, the set of $c$-clusters containing $F$ forms an interval in $\Clust_c$.
\item[(iv) ] A closed interval $I$ in $\Clust_c$ has proper part homotopy equivalent to an $(n-k-2)$-dimensional sphere if and only if there is some set $F$ of~$k$ $c$-compatible almost positive roots such that $I$ is the set of all $c$-clusters containing $F$.
\item[(v) ] A closed interval $I$ in $\Clust_c$ has proper part homotopy equivalent to an $(k-2)$-dimensional sphere if and only if $I$ has~$k$ atoms and the join of the atoms of $I$ is the top element of $I$.
\item[(vi) ] If the proper part of a closed interval $I$ is not homotopy-spherical then it is contractible.
\end{enumerate}
\end{cor}

Theorem ~\ref{cl poset iso} also has important enumerative consequences.
Let $C \subset \Pge$ be a $c$-cluster. 
For each $\alpha \in C$, there is a unique $\alpha' \in \Pge$ such that $(C \setminus \set{\alpha}) \cup \set{\alpha'}$ is also a $c$-cluster. 
Call~$\alpha$ an {\em upper root} of~$C$ if $R_c(\alpha)>R_c(\alpha')$ and a {\em lower root} of~$C$ if $R_c(\alpha) < R_c(\alpha')$.
Equivalently,~$\alpha$ is an upper root if $C\covers (C \setminus \set{\alpha}) \cup \set{\alpha'}$ or a lower root if $C\covered (C \setminus \set{\alpha}) \cup \set{\alpha'}$.
Note that the identification of~$\alpha$ as a lower or upper root depends on the $c$-cluster~$C$.  
A root may be an upper root in one $c$-cluster and a lower root in another $c$-cluster.

\begin{cor}
\label{cl h}
The map $\cl_c$ takes $c$-sortable elements with~$k$ descents to $c$-clusters with~$k$ upper roots.
\end{cor}
Since any linear extension of $\Clust_c$ is a shelling order, by standard arguments the number of $c$-clusters with~$k$ upper roots is the $k^{\textrm{th}}$ entry in the $h$-vector of $c$-cluster fan, or equivalently the $h$-vector of the generalized associahedron for~$W\!$.
Thus Corollary~\ref{cl h} combines with the second sentence of Theorem~\ref{nc} to give a bijective proof of the following.

\begin{cor}
\label{nc cl h}
The number of noncrossing partitions (with respect to~$c$) of rank~$k$ equals the $k^{\textrm{th}}$ entry in the $h$-vector of the generalized associahedron for~$W\!$.
\end{cor}
This number is known as the $k^{\textrm{th}}$ Narayana number associated to~$W\!$.
Corollary~\ref{cl h} can be viewed as a direct combinatorial interpretation of the Narayana numbers in terms of $c$-clusters.

We now proceed to prove Propositions~\ref{BigR} and~\ref{cl covers}.
We begin by proving a strengthening of Proposition~\ref{BigR}. This argument follows a suggestion of a referee. 

\begin{prop} \label{LittleR}
For any almost positive root~$\alpha$ and any reduced word $s_1\cdots s_n$ for a Coxeter element~$c$, there exists a nonnegative integer~$r$ with the property that $\sigma_{s_r} \sigma_{s_{r-1}} \cdots \sigma_{s_2} \sigma_{s_1} \alpha$ is a negative simple root.
\end{prop}
Here, the subscripts are interpreted cyclically, so that $s_{n+1}=s_1$, etc.
\begin{proof}
We first claim that the $c$-orbit of $\alpha$ contains a negative root.
Suppose to the contrary that every root in the $c$-orbit of $\alpha$ is positive.
Then in particular the sum of the roots in the $c$-orbit is a nonzero vector is fixed by $c$.  
(The assumption that $W$ is finite is critical here.  If $W$ is infinite then the $c$-orbit of $\alpha$ may be infinite, so it may not have a well-defined sum.)
However, it is well-known that $c$ acts without fixed points (see for example \cite[Section~V.6.2]{Bourbaki} or \cite[Lemma~3.16]{Humphreys}).
This contradiction proves the claim.

If $\alpha$ is a negative simple root then take $r=0$.
Otherwise let $R$ be the smallest nonnegative integer such that $c^{-R}(\alpha)$ is a negative root and let $\beta$ be the positive root $c^{-R+1}(\alpha)$.
We claim that for $r$ between 0 and $(-R+1)n$, the root $s_rs_{r-1}\cdots s_1\alpha$ is positive.
To prove the claim by contradiction, take $r$ to be the smallest exception and let $r'$ be the smallest multiple of $n$ greater than $r$.
Since each \mbox{$s\in S$} only changes the positive/negative status of the roots $\pm\alpha_s$, necessarily $s_rs_{r-1}\cdots s_1\alpha=-\alpha_{s_r}$.
Furthermore $s_{r'}s_{r'-1}\cdots s_1\alpha=s_{r'}s_{r'-1}\cdots s_{r+1}(-\alpha_{s_r})$ is a negative root.
But $s_{r'}s_{r'-1}\cdots s_1\alpha=c^{-R'}(\alpha)$ for some $R'$ with $0<R'<R$, contradicting the choice of $R$.
This proves the claim, implying in particular that $\beta=\sigma_c^{-R+1}(\alpha)$.

Since $c\beta$ is negative and since each \mbox{$s\in S$} only changes the positive/negative status of the roots $\pm\alpha_s$, there is some $i\in[n]$ such that $s_js_{j-1}\cdots s_1\beta$ is positive for all $j<i$ and $s_{i-1}s_{i-2}\cdots s_1\beta=\alpha_{s_i}$.
Setting $r=(-R+1)n+i$ we have $\sigma_{s_r} \sigma_{s_{r-1}} \cdots \sigma_{s_2} \sigma_{s_1} \alpha=s_rs_{r-1}\cdots s_1\alpha=-\alpha_{s_i}$.
\end{proof}

Let $r_{s_1 \cdots s_n}(\alpha)$ be the smallest nonnegative integer $r$ such that $\sigma_{s_r} \sigma_{s_{r-1}} \cdots \sigma_{s_1} \alpha$ is a negative simple root.
The information given by $r_{s_1 \cdots s_n}$ is more refined than that of $R_c$ and the behavior of $r_{s_1 \cdots s_n}$ is simpler to describe. However, as the notation suggests,~$r$ depends not only on~$c$ but on a choice of a reduced word for~$c$. 
The following lemma shows how Proposition~\ref{LittleR} implies Proposition~\ref{BigR} and describes the relationship between $R_c$ and $r_{s_1 \cdots s_n}$.
Its proof follows immediately from the proof of Proposition~\ref{LittleR}.

\begin{lemma} \label{DivideByN}
For any almost positive root~$\alpha$, the integer $R_c(\alpha)$ exists and equals $\lceil r_{s_1 \cdots s_n}(\alpha) /n \rceil$.
\end{lemma}

We now proceed with the proof of Proposition~\ref{cl covers}, which states that if $B\cup\set{\alpha}$ and $B\cup\set{\alpha'}$ are clusters then \mbox{$\cl^{-1}_c(B\cup\set{\alpha}) \covered$}  \mbox{$\cl^{-1}_c(B\cup\set{\alpha'})$} if and only if $R_c(\alpha)<R_c(\alpha')$.

\begin{proof}[Proof of Proposition~\ref{cl covers}.]
We shall in fact prove that, for any reduced word $s_1\cdots s_n$ for~$c$, $\cl^{-1}_c(B\cup\set{\alpha})\covered \cl^{-1}_c(B\cup\set{\alpha'})$ if and only if $r_{s_1\cdots s_n}(\alpha)<r_{s_1\cdots s_n}(\alpha')$.
By Lemma~\ref{DivideByN} and the fact that $R_c(\alpha)$ cannot equal $R_c(\alpha')$, this implies Proposition~\ref{cl covers}.

Let $w=\cl^{-1}_c(B\cup\set{\alpha})$, let $w'=\cl^{-1}_c(B\cup\set{\alpha'})$ and let $[w]_c$ and $[w']_c$ be the corresponding maximal cones in the $c$-Cambrian fan.
By Corollary~\ref{Hasse}, either $[w]_c \covered [w']_c$ or $[w]_c \covers [w']_c$. 
Since the possibility $r_{s_1\cdots s_n}(\alpha)=r_{s_1\cdots s_n}(\alpha')$ is also ruled out, it suffices by symmetry to prove one direction of implication.
Thus we will prove that if $w\covered w'$ then $r_{s_1\cdots s_n}(\alpha)<r_{s_1\cdots s_n}(\alpha')$.
We will use Lemma~\ref{only ray} repeatedly. 
Let $r=\min( r_{s_1 \cdots s_n} (\alpha),  r_{s_1 \cdots s_n} (\alpha'))$. 
Our proof is by induction on~$r$. 

First, suppose that $r=0$, so either~$\alpha$ or $\alpha'$ is a negative simple root. 
Let~$H$ be the hyperplane separating $[w']_c$ from $[w]_c$ in~$\F_c$.
Then $\phi_c^{-1}(\alpha')$ is strictly above~$H$. The rays $\rho_{s_i}$ are not strictly above any hyperplane in $\A(W)$, so $\phi_c^{-1}(\alpha') \neq \rho_{s_i}$ and $\alpha' \neq \phi_c(\rho_{s_i})=-\alpha_{s_i}$. In other words, $\alpha'$ is not a negative simple root and $r_{s_1 \cdots s_n}(\alpha')>r_{s_1 \cdots s_n}(\alpha)=0$ as desired.

Now, we consider the situation where $r \geq 1$. 
There are three cases. 
For brevity, set $s=s_1$. 
The first case is that $s\le w$ and $s\le w'$. 
Then $Z_s(w')=sw' > sw=Z_s(w)$.
By Lemma~\ref{BigZ}, $\cl_{scs}(Z_s(w))=\sigma_s(B)\cup\set{\sigma_s\alpha}$ and $\cl_{scs}(Z_s(w'))=\sigma_s(B)\cup\set{\sigma_s\alpha'}$.
By Proposition~\ref{Hasse}, one of $[Z_s(w)]_{scs}$ and $[Z_s(w')]_{scs}$ covers the other in $W/\Theta_{scs}$, but, by the isomorphism $[1, s w_0] \isomorph [s, w_0]$, we have $Z_s(w)\le Z_s(w')$ in weak order, so $[Z_s(w)]_{scs}\covered [Z_s(w')]_{scs}$.
By induction, $r_{s_2 \cdots s_n s_1}(\sigma_{s}(\alpha)) < r_{s_2 \cdots s_n s_1}(\sigma_{s}(\alpha'))$ and thus $r_{s_1 \cdots s_n}(\alpha) < r_{s_1 \cdots s_n}(\alpha')$.

The second case, $s\not\le w$ and $s\not\le w'$, is very similar to the preceding one.
By Lemma~\ref{BigZ}, $\cl_{scs}(w)$ and $\cl_{scs}(w')$ differ only by the exchange of $\sigma_s(\alpha)$ for $\sigma_s(\alpha')$. 
By Proposition~\ref{Hasse}, $[Z_s(w)]_{scs}$ and $[Z_s(w')]_{scs}$ are a covering pair in $W/\Theta_{scs}$. 
Since $Z_s(w) = s \join w \le s \join w' = Z_s(w')$, the cover must be $[Z_s(w)]_{scs} \covered [Z_s(w')]_{scs}$. 
As in the previous case, we conclude that $r_{s_1 \cdots s_n}(\alpha) < r_{s_1 \cdots s_n}(\alpha')$.

The case $s\le w$ and $s\not\le w'$ is impossible because $w' > w$. 
So we complete the proof by considering the case $s\not\le w$ and $s\le w'$. 
Then $\phi_c^{-1}(\alpha)$ and $\phi_c^{-1}(B)$ are on or below $H_s$ and $\phi_c^{-1}(\alpha')$ and $\phi_c^{-1}(B)$ are on or above $H_s$. 
Thus all the rays in $\phi_c^{-1}(B)$ are contained in $H_s$. 
Since $\phi_c^{-1}(\alpha)$ is not in the linear span of $\phi_c^{-1}(B)$, $\phi_c^{-1}(\alpha)$ must be strictly below $H_s$. 
But then by Proposition~\ref{only ray}, $\phi_c^{-1}(\alpha)=\rho_s$ and $\alpha=-\alpha_s$. 
This is the case $r=r_{s_1 \cdots s_n}(\alpha)=0$ which we have already described.
\end{proof}

For any Coxeter element~$c$ of~$W$ and any $J\subseteq S$, let $c'$ be the Coxeter element of $W_J$ obtained by deleting the letters $S\setminus J$ from any reduced word for~$c$.
Since the $c'$-Cambrian lattice $W_J/\Theta_{c'}$ is a lower interval in $W/\Theta_c$, we have the following combinatorial fact about clusters which appears to be difficult to prove directly:

\begin{prop}
\label{R para}
For~$c$ and $c'$ as above, if~$\alpha$ and $\alpha'$ are $c$-exchangeable almost positive roots then $R_{c'}(\alpha)<R_{c'}(\alpha')$ if and only if $R_c(\alpha)<R_c(\alpha')$.  
Thus a root in a $c'$-cluster~$C$ is an upper root in~$C$ if and only if it is an upper root in the $c$-cluster $C\cup\set{-\alpha_s:s\in S\setminus J}$, and the same is true for lower roots.
\end{prop}

\begin{remark}\label{face link remark}
It is known that every face of an associahedron is combinatorially isomorphic to another associahedron.
Equivalently, the link of any cone in the cluster complex is combinatorially isomorphic to a cluster complex.
One can prove a stronger version of this result in the Cambrian setting, showing that the star of a face in the Cambrian fan is not only combinatorially a Cambrian fan, but has the polyhedral and lattice structure of a Cambrian fan as well. Specifically, for $w$ any $c$-antisortable element and $J$ a set of ascents of $w$, there is a choice of Coxeter element $\ccc(w,J,c)$ such that the following proposition holds. 

\begin{prop} \label{face link prop}
Let $C_{\Theta_c}(w,J)$ be a face of the $c$-Cambrian fan. Identify\footnote{Note that we identify a cone $\kappa$ in the $(W_J)$-Coxeter fan with a cone isomorphic to $\kappa \times \RR^{n-|J|}$ in the star of $C(w,J)$.} the star of $C(w,J)$ (in the $W$-Coxeter fan) with the star of $C(e,J)$, and hence with the $W_J$-Coxeter fan, by the map $w^{-1}$. Then $\F_{\ccc(w,J,c)}$ and the star of $C_{\Theta_c}(w,J)$ coincide as coarsenings of the $W_{J}$-Coxeter fan.
\end{prop}

%SinceJEMS: The last instance of "c(w,J,c)" changed to "\ccc(w,J,c)"

Defining $\ccc(w,J,c)$ means deciding, for each $r_1$,~$r_2 \in J$ with $r_1 r_2 \neq r_2 r_1$, whether the reflection~$r_1$ comes before~$r_2$ in every reduced word for~$\ccc(w,J,c)$ or \emph{vice versa}. In \cite[Section~3]{sortable}, a directed graph is defined on the set $T$ of reflections of $W,$ with arrows $\toname{c}$. We put $r_1$ before~$r_2$ in $\ccc(w,J,c)$  if and only if $w r_1 w^{-1} \toname{c} w r_2 w^{-1}$.

To prove Proposition~\ref{face link prop}, one first reduces to the case that $C_{\Theta_c}(w,J)$ is a ray $\rho(w,J)$. For $s$ initial in $c$, one analyzes the effect of $\zeta_s$ and $Z_s$ on the star of $\rho(w,J)$. When $w \geq s$ the star is unaltered. When $w \not\geq s$ and $w \neq 1$, the star of $\rho(w,J)$ is partly below $H_s$ and partly above.
Passing from $\F_c$ to $\F_{scs}$ has the effect of swapping the part above with the part below, as explained in Example~\ref{A3Z} and illustrated in Figure~\ref{A3Zfig}. In either case, the effect is compatible with the properties of $\toname{c}$ established in \cite[Proposition~3.1]{sortable}.  By Proposition~\ref{LittleR}, one eventually reaches a ray of the dominant chamber, where the proposition is straightforward.
\end{remark}

\section{A linear isomorphism} 
\label{lin sec}
In this section we show that for a special choice of~$c$, the $c$-Cambrian fan is linearly isomorphic to the $c$-cluster fan.
We also describe, for a special choice of~$c$, a ``twisted'' version of the $c$-cluster lattice which is induced on the $c$-cluster fan by any vector in a certain cone in $\RR^n$.

Recall that~$\Phi$ is a fixed root system for~$W\!$. For $\alpha \in \Phi$, the corresponding {\em coroot} is $\alpha\ck=\frac{2\alpha}{\br{\alpha,\alpha}}$, so that the reflection of a vector~$v$ in the hyperplane perpendicular to a root~$\alpha$ is $v-\br{v,\alpha\ck}\alpha$.
The simple coroots $\alpha_s\ck$ for \mbox{$s\in S$} are a basis for~$V$ and the {\em fundamental weights} $\omega_s$ are the dual basis vectors to the simple coroots.\footnote{We abuse terminology slightly by calling these ``weights'' even in the non-crystallographic case, where there is no ``weight lattice.''} We have 
\begin{equation}
\label{root eq}
\alpha_s=\sum_{r \in S} \br{ \alpha_s, \alpha_r^{\vee}} \omega_r.
\end{equation}
We define the Cartan matrix of~$\Phi$ to be the $n \times n$ square matrix $A$ where $A_{ij}=\br{ \alpha_{s_i}\ck, \alpha_{s_j}}$. A root system is called crystallographic if all the entries of $A$ are integers.\footnote{Our convention for $A_{ij}$ is the convention used in~\cite{ga},~\cite{ca2} and~\cite{ca4}; some references use the transpose of this choice.}

Fix a bipartition $S=S_+ \sqcup S_-$ of the Coxeter diagram for $W,$ let $c_+$ be the product of the elements of $S_+$ (which commute pairwise) and let $c_-$ be the product of the elements of $S_-$ (which likewise commute). 
The Coxeter element $c_+c_-$ is called a {\em bipartite} Coxeter element.
For each \mbox{$s\in S$}, let $\ep_s$ be $+1$ if $s\in S_+$ and let $\ep_s$ be $-1$ if $s\in S_-$.

Let~$L$ be the linear map that sends a simple root $\alpha_s$ to $-\ep_s\omega_s$.
(Cf.\ \cite[Conjecture 1.4]{cambrian}.)
The map~$L$ depends on the choice of bipartition, but we suppress this dependence in our notation.
Following the notation of~\cite{ga}, for $\ep\in\set{+,-}$ let $\tau_\ep=\prod_{s\in S_\ep}\sigma_s$, where again the order of composition is unimportant.
Thus $\tau_+\tau_-=\sigma_c$ (in the sense of Section~\ref{lattice sec}) for $c=c_+c_-$.
The main result of this section is the following:

\begin{theorem}
\label{L iso}
For $c=c_+c_-$, the map~$L$ is a linear isomorphism from the $c$-cluster fan to the $c$-Cambrian fan.
As a map on rays, the map~$L$ coincides with $\phi^{-1}_c\circ\,\tau_-$.
\end{theorem}

% \begin{remark}
% \label{g vector}
% This remark concerns cluster algebras of finite type. As we will discuss further in the next section, almost positive roots correspond to {\em cluster variables} in a cluster algebra of finite type. Comparing with \cite[Proposition~11.3]{ca4}, we see that $\phi_c^{-1}$ is essentially the map taking the {\em denominator vector} of a cluster variable to the {\em $\mathbf{g}$-vector} of the same cluster variable. (The difference is that the map $E\circ\,\tau_-$ of~\cite{ca4} expresses the $\mathbf{g}$-vector in the basis of simple roots, rather than in the basis of fundamental weights.)
% Thus the following statement is true:
% For bipartite~$c$, identifying the rays of the $c$-Cambrian fan with the cluster variables via the map~$\phi_c$, the $\mathbf{g}$-vector of any cluster variable (with respect to a {\em bipartite seed}) is given by the coordinates, in the basis $\set{\omega_s}$, of a certain vector in the corresponding ray. (Namely, the fundamental vector which we will define below.)  In the next section, we will prove an analogue of this result without any bipartite hypotheses, but we will need to assume a certain conjecture (Conjecture 7.12) from~\cite{ca4}. It is interesting that, in this one case, we can prove the same result without assuming this conjecture.
% \end{remark}

We begin the proof of Theorem~\ref{L iso} with a simple lemma.  

\begin{lemma}
\label{cL}
The linear maps~$c$, $c_+$, $c_-$ and~$L$ on~$V$ satisfy the following equalities.
\begin{enumerate}
\item[(i) ]$c_+L=-Lc_-$
\item[(ii) ]$c_-L=-Lc_+$
\item[(iii) ]$c^{-1}L=Lc$
\end{enumerate}
\end{lemma}
\begin{proof}
We prove equality (i) by evaluating each side of the equality on the basis elements $\alpha_s$.
If $s\in S_-$ then $c_+L\alpha_s=c_+\omega_s=\omega_s$, with the latter equality holding because $\omega_s$ is orthogonal to $\alpha_r$ for each $r\neq s$.
On the other hand, $-Lc_-\alpha_s=-L(-\alpha_s)=\omega_s$.

If $s\in S_+$ then $c_+L\alpha_s=c_+(-\omega_s)=-s\omega_s=-\omega_s+\alpha_s$.
On the other hand 
\[-Lc_-\alpha_s=-L\left(\alpha_s-\sum_{r\in S_-}\br{\alpha_s,\alpha_r\ck}\alpha_r\right)=\omega_s+\sum_{r\in S_-}\br{\alpha_s,\alpha_r\ck}\omega_r.\]
To see that these two sides are equal, we must show that
$$\alpha_s=2 \omega_s + \sum_{r \in S_{-}} \br{ \alpha_s, \alpha_r\ck} \omega_r.$$
We have $\br{ \alpha_s, \alpha_s\ck}=2$ and $\br{\alpha_s, \alpha_r\ck}=0$ for $r \in S_+ \setminus \{ s\}$, so the right hand side is $\sum_{r \in S} \br{ \alpha_s, \alpha_r\ck} \omega_r$, which, as already noted, equals $\alpha_s$. 

Now (ii) follows by reversing the roles of ``$+$'' and ``$-$'' and (iii) follows by combining (i) and (ii), keeping in mind that $c^{-1}=c_-c_+$.
\end{proof}

We now prove a version of equality (iii) in the previous lemma which is more complicated in the sense that it involves maps which are not linear.
Specifically, it uses the maps~$\zeta_s$ and~$\sigma_s$ which appear in Lemma~\ref{zeta}.
In what follows, we apply the maps~$\zeta_s$ to vectors rather than rays.
To do this, we define the {\em fundamental vector} in a ray $\rho$ of the Coxeter fan to be the unique $\omega$ in the~$W$-orbit of the fundamental weights $\set{\omega_s:s\in S}$ such that $\omega\in\rho$.
Notice that Lemma~\ref{zeta} applies even when rays are replaced by fundamental vectors.
For $\ep\in\set{+,-}$ let $\zeta_{c_\ep}=\prod_{s\in S_\ep}\zeta_s$.
For $c=c_+c_-$, let $\zeta_c=\zeta_{c_+}\zeta_{c_-}$ and $\zeta_{c^{-1}}=\zeta_{c_-}\zeta_{c_+}$.

\begin{lemma}
\label{zetaL}
For $c=c_+c_-$, if~$\alpha$ is a positive root then $L\sigma_c\alpha=\zeta_{c^{-1}}L\alpha$.
\end{lemma}
\begin{proof}
If~$\alpha$ is a simple root $\alpha_s$ for $s\in S_-$ then 
\[L\sigma_c\alpha=L\tau_+\tau_-\alpha=L\tau_+(-\alpha_s)=L(-\alpha_s)=-\omega_s.\]
On the other hand, 
\[\zeta_{c^{-1}}L\alpha=\zeta_{c_-}\zeta_{c_+}\omega_s=\zeta_{c_-}\omega_s=-\omega_s.\]

If~$\alpha$ is a positive root not of the form $\alpha_s$ for $s\in S_-$ then 
\[L\sigma_c\alpha=L\tau_+\tau_-\alpha=L\tau_+c_-\alpha=Lc_+c_-\alpha.\]
(The second equality holds because~$\alpha$ is a positive root. 
The only positive roots which are sent to negative roots by $c_-$ are roots of the form $\alpha_s$ for $s\in S_-$.
Thus $c_-\alpha$ is a positive root and therefore the third equality holds as well.)
On the other hand, since $L^{-1}(\omega_s)=-\alpha_s$ when $s\in s_+$, the vector $L\alpha$ is not of the form $\omega_s$ for $s\in S_+$.  Thus $\zeta_{c^{-1}}L\alpha=\zeta_{c_-}\zeta_{c_+}L\alpha=\zeta_{c_-}c_+L\alpha$.
If $c_+L\alpha$ is $\omega_s$ for some $s\in S_-$ then $L\alpha=\omega_s$ as well, so that $\alpha=\alpha_s$.  
Since we are currently in the case which excludes such an~$\alpha$, we can write $\zeta_{c_-}c_+L\alpha=c_-c_+L\alpha$.
Thus in this case the requirement is that $Lc\alpha=c^{-1}L\alpha$, which was proved in Lemma~\ref{cL}.
\end{proof}

The map $\phi_c^{-1}$, as defined in Section~\ref{ray sec}, takes almost positive roots to rays.
In what follows, we continue to identify each ray $\rho$ with the fundamental vector in $\rho$. 

\begin{prop}
\label{linear}
For $c=c_+c_-$, the map~$L$ takes almost positive roots to rays of the $c$-Cambrian fan.
Specifically,~$L$ restricted to almost positive roots is $\phi^{-1}_c\circ\,\tau_-$.
\end{prop}
\begin{proof}
Let~$\alpha$ be an almost positive root.
We show by induction on $R_{c^{-1}}(\alpha)$ that $L\alpha$ is a ray of the $c$-Cambrian fan and that $\phi_cL\alpha=\tau_-\alpha$.
First suppose that $R_{c^{-1}}(\alpha)=0$, so that~$\alpha$ is a negative simple root $-\alpha_s$.
In this case, $L\alpha=\pm\omega_s$, which in either case is a ray of the $c$-Cambrian fan. 
If $s\in S_+$ then $\tau_-\alpha=\alpha$ and $\phi_cL\alpha=\phi_c\omega_s=-\alpha_s=\alpha$.
If $s\in S_-$ then 
\[\phi_cL\alpha=\phi_c(-\omega_s)=\phi_c(\zeta_{c_-}\omega_s)=\tau_-\phi_{c^{-1}}\omega_s=\tau_-(-\alpha_s).\]
Here the next-to-last equality follows from Lemma~\ref{zeta}, applied several times, and the fact that $c^{-1}=c_-c_+$.

Next suppose that  $R_{c^{-1}}(\alpha)>0$ so that~$\alpha$ is a positive root and $\sigma_c\alpha=\alpha'$ for some $\alpha'$ with $R_{c^{-1}}(\alpha')=R_{c^{-1}}(\alpha)-1$.
By induction, $L\alpha'$ is a ray of the $c$-Cambrian fan and $\phi_cL\alpha'=\tau_-\alpha'$.
To evaluate $\phi_cL\alpha$, first note that by Lemma~\ref{zetaL}, 
\[L\alpha=\zeta_{c^{-1}}^{-1}L\sigma_c\alpha=\zeta_{c^{-1}}^{-1}L\alpha'.\]
In particular, $L\alpha$ is a ray in the $c$-Cambrian fan and $\phi_cL\alpha=\phi_c\zeta_{c^{-1}}^{-1}L\alpha'$.
Repeated applications of Lemma~\ref{zeta} give the identity $\phi_c\zeta_{c^{-1}}=\sigma_{c^{-1}}\phi_c$, so that $\phi_c\zeta_{c^{-1}}^{-1}=\sigma_{c^{-1}}^{-1}\phi_c=\sigma_c\phi_c$.
Thus 
\[\phi_cL\alpha=\sigma_c\phi_cL\alpha'=\sigma_c\tau_-\alpha'=\tau_+\alpha'=\tau_+\sigma_c\alpha=\tau_-\alpha.\]
\end{proof}

The map $\tau_-$ induces a combinatorial isomorphism between the $(c_+c_-)$-cluster fan and the $(c_-c_+)$-cluster fan.
These fans in fact coincide, so that $\tau_-$ is a combinatorial automorphism of the $(c_+c_-)$-cluster fan.
By Theorem~\ref{cl iso}, $\phi^{-1}_{c_+c_-}$ induces a combinatorial isomorphism as well.
This completes the proof of Theorem~\ref{L iso}.

We conclude the section with an application of Theorem~\ref{L iso}.
Proposition~\ref{linear} suggests the definition of a ``twisted'' cluster lattice on $c$-clusters, where $c=c_+ c_-$.
Namely, for $c$-clusters~$C$ and~$C'$, set $C\letw C'$ in the twisted $c$-cluster lattice if and only if $\tau_-C\le\tau_-C'$ in the $c$-cluster lattice.
In particular, the cover relations in the twisted $c$-cluster lattice are $B\cup\set{\alpha}\covered B\cup\set{\alpha'}$ if and only if $R_c(\tau_-\alpha)<R_c(\tau_-\alpha')$.

The twisted $c$-cluster lattice can be described in terms of a quantity $\ep(\alpha,\alpha')$ which plays an important role in~\cite{ca2}, where cluster algebras of finite type are constructed in terms of the combinatorics of clusters of almost positive roots.
Let $\tau_-^{(k)}$ denote the $k$-fold composition $\tau_{(-1)^k}\tau_{(-1)^{k-1}}\cdots\tau_-\tau_+\tau_-$.
For each almost positive root~$\alpha$, let $k_-(\alpha)$ be the smallest nonnegative integer such that $\tau_-^{(k)}(\alpha)$ is a negative simple root and $\tau_-^{(k)}(\alpha)=\tau_-^{(k+1)}(\alpha)$.
Given two $c$-clusters $B\cup\set{\alpha}$ and $B\cup\set{\alpha'}$, define $\ep(\alpha,\alpha')$ to be $-1$ if $k_-(\alpha)<k_-(\alpha')$ or~$1$ if $k_-(\alpha')<k_-(\alpha)$.
(See \cite[Lemma~4.1]{ca2}.)
As with $R_c$ and $r_{s_1\cdots s_n}$, the case $k_-(\alpha)=k_-(\alpha')$ is impossible.

The following proposition says that the twisted $c$-cluster lattice is analogous to the ordinary $c$-cluster lattice, except that $k_-$ plays the role of $R_{c_+c_-}$.
The proof is a straightforward induction on $k_-(\alpha)$, and we omit the details.

\begin{prop}
\label{k ep}
For $c=c_+c_-$, if~$\alpha$ and $\alpha'$ are $c$-exchangeable then
\[R_c(\tau_-\alpha)<R_c(\tau_-\alpha')\quad\mbox{if and only if }\quad\ep(\alpha,\alpha')=-1.\]
In particular, the cover relations of the twisted $c$-cluster lattice are of the form \mbox{$B\cup\set{\alpha}\covered B\cup\set{\alpha'}$} for $\ep(\alpha,\alpha')=-1$.
\end{prop}
Since the twisted $c$-cluster lattice is isomorphic to the ordinary $(c_-c_+)$-cluster lattice by a map which also induces a combinatorial isomorphism of fans, the twisted $c$-cluster lattice inherits all of the properties listed in Corollary~\ref{consequences}.
(These properties are all combinatorial.)
Since the isomorphism between the twisted $c$-cluster lattice and the $c$-Cambrian lattice is given by a linear map of fans, the following property of the $c$-Cambrian lattice (see Section~\ref{camb fan sec}) carries over to the twisted $c$-cluster lattice:
\begin{prop}
\label{twist fun}
The twisted $c$-cluster lattice is the order induced on the maximal cones of the $c$-cluster fan by any vector in the interior of the cone spanned by the $c$-cluster $\set{-\ep_s\alpha_s:s\in S}$.
\end{prop}

\section{Connections to cluster algebras} \label{gVectorSection}
In this section we connect our results to the theory of cluster algebras. 
Rather than give the lengthy definition of a cluster algebra, we merely describe the properties of cluster algebras and refer the reader to~\cite{ca4} for definitions.

Let $\FFF$ be a field isomorphic to $\QQ(x_1, \ldots, x_n)$ and let $B$ be an $n \times n$ integer matrix that is skew-symmetrizable, meaning that there exists an invertible diagonal matrix $D$ such that $DB$ is skew-symmetric. 
The combinatorial data for a cluster algebra is the matrix $B$ and an $n$-tuple $(x_1, \ldots, x_n)$ of rational functions generating $\FFF$  as a field.
The cluster algebra $\Alg(B, (x_1, \ldots, x_n))$ is a certain subalgebra of the Laurent-polynomial ring $\ZZ[x_1^{\pm}, \ldots, x_n^{\pm}]$ which is, in turn, a subring of $\FFF$. 
The data $(B, (x_1, \ldots, x_n))$ also determines a collection of transcendence bases of $\Alg(B, (x_1, \ldots, x_n))$, known as algebraic clusters.\footnote{Typically, these are simply called ``clusters,'' but we use the adjective ``algebraic'' here to avoid confusion with $c$-clusters.} 
The elements of the algebraic clusters are known as cluster variables.
One algebraic cluster is $(x_1, \ldots, x_n)$ and the others are defined by a certain recursive procedure. 
The recursive procedure also associates a skew-symmetrizable matrix $B^t$ to each algebraic cluster $t=(y_1, \ldots, y_n)$ so that $\Alg(B, (x_1, \ldots, x_n))=\Alg(B^t, (y_1, \ldots, y_n))$ and so that $(B, (x_1, \ldots, x_n))$ and $(B^t,(y_1, \ldots, y_n))$ each give the same collection of algebraic clusters.

A cluster algebra is of finite type if it has finitely many cluster variables.
We now briefly describe the connection between cluster algebras of finite type and finite Coxeter groups/root systems.
For more details, see~\cite{ca2}.
Let~$\Phi$ be a crystallographic root system for the Coxeter group~$W\!$. 
We refer the reader to the beginning of Section~\ref{lin sec} for our conventions regarding roots, coroots, Cartan matrices and fundamental weights.
Let~$c$ be a Coxeter element of~$W\!$. 
If~$r$ and~$s$ are two simple reflections of~$W$ which do not commute, then either~$r$ comes before~$s$ in every reduced word for~$c$ or \emph{vice versa}. 
We write $r \to s$ to indicate that~$r$ comes before~$s$ in every reduced word for~$c$. 

Define a square matrix $B^{c}$ by 
$$B^c_{jk} = \left\lbrace\begin{array}{rcrl}
0&&& \mbox{if } s_j s_k=s_k s_j \\[1 mm]
- A_{jk}&=&-\br{\alpha_{s_j} \ck, \alpha_{s_k}} &\mbox{if } s_j \to s_k \\[1 mm] 
A_{jk}&=&\br{\alpha_{s_j} \ck, \alpha_{s_k}} &\mbox{if } s_j \leftarrow s_k. \\ 
\end{array}\right. $$
Then the matrix $B^c$ (together with any choice of $(x_1, \ldots, x_n)$) defines a cluster algebra of finite type.
Furthermore, cluster algebras arising from different choices of $c$ and $(x_1, \ldots, x_n)$ are isomorphic; we thus suppress the choice of $c$ and $(x_1, \ldots, x_n)$ and write $\Alg(\Phi)$ for a cluster algebra arising in this manner.

Conversely, given any cluster algebra of finite type, there exists\footnote{Most often, $c$ does not uniquely determine $t_c$. Here we assume that some choice of $t_c$ has been made.  On the other hand, outside of rank two, not every algebraic cluster can serve as $t_c$.} a finite Coxeter group $W$ (with root system $\Phi$), a Coxeter element $c$ in $W$ and an algebraic cluster $t_c=(x_1^c, \ldots, x_n^c)$ such that the given cluster algebra is $\Alg(\Phi)=\Alg(B^c, (x_1^c, \ldots, x_n^c))$.
Thus the cluster algebras of finite type are precisely the cluster algebras of the form $\Alg(\Phi)$, so that the following theorem applies to any cluster algebra of finite type.

\begin{theorem} \label{cluster summary}
Given a specific representation of $\Alg(\Phi)$ as $\Alg(B^c, t_c)$, there is a bijection $\alpha \mapsto x^c(\alpha)$ between $\Phi_{\geq -1}$ and the cluster variables of $\Alg(\Phi)$ such that:
\begin{enumerate}
\item[(i) ]$x^c(- \alpha_{s_j})=x^c_j$ for all $j\in[n]$;
\item[(ii) ]$c$-clusters are mapped to algebraic clusters; and
\item[(iii) ]For positive roots $\alpha=\sum a_i \alpha_{s_i}$, the rational function $x^c(\alpha)$ can be written in reduced form with denominator $\prod x^c(-\alpha_{s_i})^{-a_i}$.
\end{enumerate}
Furthermore, if $s$ is initial in $c$ then $t_{scs}$ can be chosen so that $x^c(\alpha)=x^{scs}(\sigma_s(\alpha))$.
\end{theorem}

\begin{proof}
In the case of bipartite $c$, the first assertion is \cite[Theorem~1.9]{ca2}. 
For general~$c$, the entire theorem was proven for simply laced root systems (\emph{i.e.} $A_{ij}=0$ or $-1$ for all $i \neq j$) in~\cite{ccs}, relying on previous work cited therein. 
The result for non-simply laced root systems can be established by folding arguments. 
\end{proof}

In rough terms, Theorem~\ref{cluster summary} says that the cluster variables of $\Alg(\Phi)$ correspond to almost positive roots by assigning a variable to its \emph{denominator vector} $\prod x^c(-\alpha_{s_i})^{-a_i}$. 
There is another natural way to encode cluster variables by integer vectors, namely the $\mathbf{g}$-vector, defined in~\cite{ca4}. 
The $\mathbf{g}$-vector of a cluster variable~$x$ depends on a fixed algebraic cluster~$t$ and is written $g^{t}(x)$, with components $g_j^{t}(x)$. In~\cite[Proposition~11.3]{ca4}, Fomin and Zelevinsky compute the $ \mathbf{g}$-vector when (in the language of the current paper) $t$ 
is of the form $t_c$ for $c$ a bipartite Coxeter element. 
They encode the  $\mathbf{g}$-vector as an element of $V$ by the sum $\groot^{t}(x):=\sum_{j=1}^n g_j^{t}(x) \alpha_{s_j}$. 
(In~\cite{ca4}, this sum is also denoted by $g^{t}(x)$.
However, it is important here to distinguish between the integer vector $g^t(x)$ and the vector $\groot^{t}(x)$ lying in the root lattice.)
They establish the formula
$$\groot^{t_c}(x^c(\alpha))=(E \circ \tau_{-})(\alpha).$$
Here $\tau_{-}$ has the same meaning as it did in section~\ref{lin sec} and $E$ is the linear map such that $E(\alpha_s)=- \epsilon(s) \alpha_s$. 

The $\mathbf{g}$-vector has no obvious connection to the geometry of the $c$-cluster fan, but remarkably, it arises naturally in the geometry of the $c$-Cambrian fan.
To see this, we encode the $\mathbf{g}$-vector in the weight lattice by $\gweight^t(x):=\sum_{j=1}^n g_j^{t}(x) \omega_{s_j}$.
Let $U$ denote the linear map which takes $\alpha_s$ to $\omega_s$, so that $\gweight^t(x)= U(\groot^t(x))$. 
Thus when $c$ is bipartite, Theorem~\ref{L iso} implies that
\[\gweight^{t_c}(x^c(\alpha))=(U \circ E \circ \tau_{-})(\alpha)=(L \circ \tau_{-})(\alpha)=\phi_c^{-1}(\alpha).\]

\begin{theorem} \label{bipartite g}
If $c$ is a bipartite Coxeter element, with $t_c$ a corresponding cluster, and if $\alpha$ is an almost positive root, then
 $$\phi_c^{-1}(\alpha)=\gweight^{t_c}(x^c(\alpha)).$$
\end{theorem}
Thus $\mathbf{g}$-vectors arise naturally from the correspondence between cluster variables and rays in the Cambrian fan:
the $\mathbf{g}$-vector associated to a ray is recovered by computing the fundamental-weight coordinates of the fundamental vector in the ray.
This is precisely analogous to the situation in the cluster fan, where the denominator vector associated to a ray is recovered by taking the simple-root coordinates of the root in the ray.

We conjecture that Theorem~\ref{bipartite g} is true without assuming that $c$ is bipartite. 
In \cite[Conjecture 7.12]{ca4}, Fomin and Zelevinsky give a conjectured recurrence for $g^{t}(x)$ as $t$ varies. 
By a straightforward but lengthy computation, one can verify that the more general version of Theorem~\ref{bipartite g} follows from \cite[Conjecture 7.12]{ca4}.

We now sketch an additional connection between Cambrian fans and cluster algebras.
Choose a Coxeter element~$c$ of $W$ and a cluster $t_c$ of $\Alg(\Phi)$ as above.
Let $t=(x_1,\ldots, x_n)$ be an arbitrary algebraic cluster in $\Alg(\Phi)$.
(In particular, we do \emph{not} assume that~$B^t=B^{c'}$ for some $c'$.)
Then  $x_i=x^c(\alpha_i)$ for some $c$-cluster $(\alpha_1, \ldots, \alpha_n)$.
Let $[w]_c$ be the cone of $\F_c$ represented by a $c$-sortable element~$w$ with $\cl_c(w)=(\alpha_1, \ldots, \alpha_n)$.
There is another collection of roots, besides the $\alpha_i$, naturally associated to $[w]_c$, namely the roots $(\beta_1, \ldots, \beta_n)$ orthogonal to the walls of $[w]_c$.
More specifically, let $\beta_i$ be the root determined by the requirements that  $\br{\phi_c^{-1}(\alpha_i), \beta_j }=0$ for $i \neq j$ and $\br{\phi_c^{-1}(\alpha_i), \beta_i} < 0$. Let~$\beta_i\ck$ be the coroot corresponding to the root $\beta_i$ and let $Q^t$ be the $n \times n$ matrix $\br{ \beta_i\ck, \beta_j}$.

The matrix $Q^t$ depends on the choice of $c$ and $t_c$ above. 
This dependence is not as bad as one might suspect. 
If $s$ is initial in $c$ and $t_{scs}$ is the cluster referred to in Theorem~\ref{cluster summary}, then the cone corresponding to $t$ changes from $[w]_c$ to $[Z_s(w)]_{scs}$. 
Either $[w]_c$ and $[Z_s(w)]_{scs}$ are related by an isometry or else they are two regions among the $2^n$ regions defined by the same set of $n$ hyperplanes. 
In the first case, $Q^t$ is preserved, in the second it is conjugated by a diagonal matrix all of whose diagonal entries are~$\pm 1$. 
So changing $t_c$ in this manner any number of times simply conjugates $Q^t$ by such a matrix. 
We will see soon that it follows from results of \cite{capsm} that $Q^t$ is well defined up to such conjugation independent of any of our choices.\footnote{In a previous version of this paper, we argued that this independence could be established by a sequence of steps, each changing from $t_c$ to $t_{scs}$ for $s$ initial.  A comment by one of the referees has lead us to doubt this argument.  It follows from~\cite[Theorem~1.2(1)]{BGP} that we may change $c$ to any other Coxeter element $c'$ by such a sequence of steps.
What is not clear is whether we may change any cluster $t_c$ corresponding to $c$ to any cluster $t_{c'}$ corresponding to $c'$.}

\begin{prop}
\label{quasi}
Let $Q^t$ be as above and let $B^t$ be the matrix associated to the algebraic cluster~$t$.
Then $Q^t_{ij} = \pm B_{ij}^t$ for $i \neq j$.
\end{prop}

\begin{proof}[Sketch of Proof:]
For $i \neq j$, the quantity $B_{ij}^t B_{ji}^t$ is encoded in the combinatorics of the algebraic cluster complex:
One counts the number of algebraic clusters containing $t \setminus \{ x^t_i, x^t_j \}$.
This number is $4$, $5$, $6$ or $8$, corresponding (in order) to $B_{ij}^t B_{ji}^t=0$, $-1$, $-2$  or $-3$.
To prove Proposition~\ref{quasi}, we verify that $-Q^t_{ij} Q^{t}_{ji}$ takes only the values $0$, $-1$, $-2$  or $-3$ and that the value of $-Q^t_{ij} Q^{t}_{ji}$ corresponds to the number ($4$, $5$, $6$ or $8$) of $c$-clusters containing $(\alpha_1, \ldots, \alpha_n) \setminus \{ \alpha_i, \alpha_j \}$.
Once this is verified, we have $B_{ij}^t B_{ji}^t = - Q_{ij}^t Q_{ji}^t$ by the isomorphism between the $c$-cluster complex and the algebraic cluster complex. The matrices $Q^t$ and $B^t$ are (respectively) symmetrizable and skew-symmetrizable. One can check that the same diagonal matrix $D$ makes both $D Q^t$ symmetric and $D B^t$ skew symmetric so we conclude from $B_{ij}^t B_{ji}^t = - Q_{ij}^t Q_{ji}^t$ that $Q^t_{ij} = \pm B_{ij}^t$.

Let $F$ be the face of $[w]_c$ spanned by $\phi_c^{-1}\left( (\alpha_1, \ldots, \alpha_n) \setminus \{ \alpha_i, \alpha_j \} \right)$. In other words, $F$ is $[w]_c \cap \beta_i^{\perp} \cap \beta_j^{\perp}$. 
By Theorem~\ref{cl iso}, the number of $c$-clusters containing $ (\alpha_1, \ldots, \alpha_n) \setminus \{ \alpha_i, \alpha_j \}$ is equal to the number of maximal faces of the $c$-Cambrian fan containing $F$.
By Proposition~\ref{face link prop}, the star of $F$ is a Cambrian fan for a (crystallographic) Coxeter group of rank $2$, of which there are only four types.
Moreover, $\beta_i$ and $\beta_j$ are roots in a rank $2$ root subsystem of corresponding type. 
(Specifically, if $F=C_{\Theta_c}(w,J)$ then $\beta_i$ and $\beta_j \in w \Phi_{J}$.) 
By inspection of Cambrian lattices of rank $2$, we see that either $(\beta_i, \beta_j)$ or $(\beta_i, - \beta_j)$ form a simple system for this root subsystem.
Thus $- \br{\beta_i^{\ck}, \beta_j} \br{\beta_j^{\ck}, \beta_i}$ is $0$, $-1$, $-2$ or $-3$ according to whether the root subsystem is $A_1 \times A_1$, $A_2$, $B_2$ or $G_2$; this in turn corresponds to whether the star of $F$ has $4$, $5$, $6$ or $8$ maximal cones.
 \end{proof}

Rephrased in the language of~\cite{capsm}, Proposition~\ref{quasi} says that $Q^t$ is a {\em quasi-Cartan companion} for $B^t$. 
The matrix $Q^t$ is positive definite\footnote{More accurately, $D Q^{t}$ is positive definite, but we follow the convention of~\cite{capsm} of saying that $Q^t$ is positive definite in this case.} because it is (essentially) a matrix of inner products between $n$ linearly independent vectors.
One direction of \cite[Theorem~1.2]{capsm} states that $B^t$ has a positive definite quasi-Cartan companion which is (by \cite[Propositions~1.4 and~1.5]{capsm}) unique up to conjugation by diagonal matrices with diagonal entries $\pm 1$.
Thus $Q^t$ is unique up to such conjugation. 
One of the virtues of this manner of obtaining $Q^t$ is that this uniqueness occurs for a geometrically natural reason, as described above. 

In this section we have suggested two new geometric approaches to the study of cluster algebras. 
First, to encode $\mathbf{g}$-vectors as linear combinations of the fundamental weights. Second, to view quasi-Cartan companions as matrices of inner products between normal vectors to a simplicial cone in the hyperplane arrangement. We hope that both of these ideas will have wider applications, including applications beyond finite type.

\section{Clusters and noncrossing partitions}
\label{bij sec}
In light of Theorems~\ref{nc} and~\ref{cl}, the map $\nc_c\circ\cl_c^{-1}$ is a bijection from $c$-clusters to $c$-noncrossing partitions.
In this section we describe this composition as a direct map, eliminating the intermediate $c$-sortable elements.
For brevity, we continue to leave out the precise details about noncrossing partitions. The $c$-noncrossing partitions are certain elements of $W$. Brady and Watt showed~\cite [Lemma 5]{BWorth} that a $c$-noncrossing partition can be recovered (among the set of all $c$-noncrossing partitions) from its fixed point set. The fixed point set of a $c$-noncrossing partition is called a \emph{$c$-noncrossing subspace}. The map $\nc_c$ takes the cover reflections of a $c$-sortable element~$w$ and multiplies them in a certain specific order such that the result is a $c$-noncrossing partition. Let $\NC_c$ be the map taking a $c$-sortable element~$w$ to the fixed points of $\nc_c(w)$; this is the intersection of the reflecting hyperplanes associated to cover reflections of~$w$.
The map $\NC_c$ is a bijection between $c$-sortable elements and $c$-noncrossing subspaces.

The composition $\NC_c\circ\cl_c^{-1}$ takes a $c$-cluster~$C$ to the intersection $I$ of the set of hyperplanes separating $[\cl_c^{-1}(C)]_c$ from equivalence classes which it covers in $W/\Theta_c$.
The subspace $I$ equals the linear span of the rays of $[\cl_c^{-1}(C)]_c$ contained in $I$.
A ray is in $I$ if and only if it is not an upper root of~$C$, i.e.\ if and only if it is a lower root of~$C$.
Thus 
\begin{theorem}
\label{yes phi}
The bijection $\NC_c\circ\cl_c^{-1}$ maps a $c$-cluster~$C$ to the $c$-noncrossing subspace $\Span_{\RR} \set{\phi^{-1}_c(\alpha):\alpha\mbox{ is a lower root in }C}$.
In particular the $c$-cluster~$C$ is uniquely identified by this subspace.
\end{theorem}

This description of the bijection has the disadvantage of depending on the recursively defined function~$\phi_c$ and on a notion of lower roots in clusters which is also defined recursively. We conjecture the following description of $I$, which would eliminate the map~$\phi_c$. 

\begin{conj} \label{span vs intersection}
Let $(\alpha_1, \ldots, \alpha_k, \beta_1, \ldots, \beta_{n-k})$ be a $c$-cluster of $W,$ with $\alpha_i$ the lower roots and $\beta_i$ the upper roots. Then 
\[\Span_{\RR} (\phi_c^{-1}(\alpha_1), \ldots, \phi_c^{-1}(\alpha_k)) = \beta_1^{\perp} \cap \cdots \cap \beta_{n-k}^{\perp}.\]
\end{conj}

It is easy to see that $\beta_1^{\perp} \cap \cdots \cap \beta_{n-k}^{\perp}$ and $\Span_{\RR} (\phi_c^{-1}(\alpha_1), \ldots, \phi_c^{-1}(\alpha_k))$ have the same dimension, so in order to prove Conjecture~\ref{span vs intersection} it is enough to show that the former contains the latter, i.e.\ that $\phi_c^{-1}(\alpha_i) 
\perp \beta_j$ for all $i$ and $j$. 
This orthogonality has been verified computationally for all choices of~$c$ in all Coxeter groups whose rank is at most~7.
Combined with Theorem~\ref{yes phi}, Conjecture~\ref{span vs intersection} would immediately imply the following conjecture.

\begin{conj} \label{no phic}
Consider the map taking a $c$-cluster~$C$ to the intersection of the hyperplanes orthogonal to the upper roots of~$C$.
This map is a bijection from $c$-clusters to $c$-noncrossing subspaces.  It coincides with $\NC_c\circ\cl_c^{-1}$.
\end{conj}

In the case of bipartite $c=c_+c_-$, the related bijection $\NC_c\circ\cl_c^{-1}\circ\,\tau_-$ can be described in a completely geometric manner, as follows. 
Recall from Section~\ref{camb fan sec} the definition of the bottom face, with respect to a generic vector, of a maximal cone in a simplicial fan.
Choosing a vector $v$ as in Proposition~\ref{twist fun}, we map each maximal cone~$C$ to the subspace $\Span_\RR(L(F))$, where~$L$ is the linear map of Section~\ref{lin sec} and $F$ is the the bottom face of~$C$ with respect to~$v$. 
In light of Propositions~\ref{linear} and~\ref{twist fun}, $\Span_\RR(L(F))$ is the span of $\set{\phi_c^{-1}\tau_-\alpha:\tau_-\alpha\mbox{ is a lower root in }\tau_-C}$.
Thus 
\begin{theorem}\label{geom bij}
The map $C\mapsto \Span_\RR(L(F))$ is the bijection $\NC_c\circ\cl_c^{-1}\circ\,\tau_-$ from $(c_+c_-)$-clusters to $(c_+c_-)$-noncrossing subspaces. 
In particular a maximal cone~$C$ in the $(c_+c_-)$-cluster complex is uniquely determined by $\Span_\RR(F)$.
\end{theorem}

\begin{remark}
\label{ABMW remark}
In~\cite{ABMW}, Athanasiadis, Brady, McCammond and Watt give another bijection between $c$-clusters and $c$-noncrossing partitions for the case where~$c$ is bipartite.
A key element of their bijection is a labeling of the roots of each cluster as ``left'' or ``right'' roots \cite[Section~4]{ABMW}.
The cluster is then mapped to the product of the reflections corresponding to its right roots, in some specified order.
Although the connection is not immediately obvious, it is natural to suspect that the left-right dichotomy of~\cite{ABMW} corresponds to the upper-lower dichotomy of the present paper.
In particular, it seems quite likely that the map of~\cite{ABMW} coincides, in the bipartite case, with the bijection of Conjecture~\ref{no phic}.
\end{remark}

\section*{Acknowledgments}

We are grateful to Sergei Fomin, John Stembridge and Andrei Zelevinsky for many helpful conversations. We would also like to thank the American Institute of Mathematics, at which our collaboration began, for their excellent hospitality. 
Finally, we thank the anonymous referees for suggestions which improved the exposition of the paper and for a mathematical suggestion which led to a significant simplification of the proofs in Section~\ref{lattice sec}.

\end{document}